%% file: partact.tex
\documentclass[twoside]{elsart}

\def\MSCnew{\par\leavevmode\hbox {\it 2000 MSC:\ }}%

\usepackage{amsmath,amssymb}
\input{header}

\input{diagrams}

\diagramstyle[size=1cm]
\newarrow{Into} C--->
\newarrow{Line} -----
\newarrow{Myto} ---->
\newarrow{Dotsto} ....>

\newcommand{\eat}[1]{}
\newcommand{\ovX}{\overline{X}}
\newcommand{\ovd}{\overline{d}}
\newcommand{\td}{\widetilde{d}}

\begin{document}

\begin{frontmatter}

\title{Globalization of Confluent Partial Actions on Topological and Metric
Spaces}

\author[BarIlan]{Michael Megrelishvili} and
\author[Bremen]{Lutz Schr{\"o}der}
\address[BarIlan]{Dept. of Math., Bar-Ilan University, 52900 Ramat-Gan, Israel\ead{megereli@math.biu.ac.il}}
\address[Bremen]{BISS, Dept. of Comput. Sci., Bremen University, 
P.O.~Box 330 440, 28334~Bremen,
Germany\ead{lschrode@informatik.uni-bremen.de, Phone: +49-421-218-4683,
Fax: +49-421-218-3054
}}
\begin{abstract}
We generalize Exel's notion of partial group action to monoids. For 
partial monoid actions that can be defined by means of suitably
well-behaved systems of generators and relations,
we employ classical rewriting theory in order to describe the universal
induced global action on an extended set. This universal action can be
lifted to the setting of topological spaces and continuous maps, as well
as to that of metric spaces and non-expansive maps. Well-known
constructions such as Shimrat's homogeneous extension are special cases
of this construction. We investigate various properties of the arising
spaces in relation to the original space; in particular, we prove
embedding theorems and preservation properties concerning separation
axioms and dimension. These results imply that every normal (metric)
space can be embedded into a normal (metrically) ultrahomogeneous space
of the same dimension and cardinality.

\end{abstract}

\begin{keyword}
Partial action \sep ultrahomogeneous space \sep rewriting \sep 
globalization

\MSCnew 54D35 
\sep 54D15 
\sep 54F45 
\sep 20M30 
\sep 20M05 
\end{keyword}

\end{frontmatter}

\markboth{Michael Megrelishvili and Lutz Schr\"oder}{Globalization of
Confluent Partial Actions}

\section*{Introduction}

Many extension problems in topology involve the question whether a given
collection of partial maps on a space can be realized as the set of
traces of a corresponding collection of total maps on some
superspace. Consider, for example, the problem of constructing a
homogeneous extension of a given topological (or metric) space $X$. A
space is homogeneous (ultrahomogeneous) iff each partial homeomorphism
(isometry) between two singleton (finite) subspaces extends to a global
homeomorphism (isometry)~\cite{CGK} (Cf.\
also~\cite{Bogatyi02,KechrisEA04,Pestov99}, and~\cite{Cameron00} for
ultrahomogeneous graphs). One way to look at the extension problem is to
regard these partial maps as algebraic operators, so that we have a set
of generators and relations for an algebra; the algebra thus generated
can be expected to serve as a carrier set for the extended
space. Indeed, this is precisely what happens in the constructions by
Shimrat~\cite{Shimrat54}, Belnov~\cite{Belnov78},
Okromeshko and Pestov~\cite{OkromeshkoPestov84},
Uspenskij~\cite{Uspenskij98}, and
Megrelishvili~\cite{Megrelishvili85,Megrelishvili86,Megrelishvili00}.

  
Here, we pursue this concept at what may be hoped is the right level of
generality: we begin by providing a generalization of Exel's notion of
\emph{partial group action}~\cite{Exel98} to partial actions of monoids
(i.e.\ the elements of the monoid act as partial maps on the space; cf.\
Definition~\ref{def:monoidpa}). Partial actions of monoids are
characterized in the same way as partial group actions as restrictions
of global actions to arbitrary subsets. We then study properties of the
\emph{globalization} of a partial action, i.e.\ of the extended space
which is universal w.r.t.\ the property that it has a global action of
the original monoid. Most of the results we obtain depend on
\emph{confluence} of the partial action.  Here, confluence means that
the monoid and the carrier set of the globalization are given in terms
of generators and relations in such a way that equality of elements can
be decided by repeated uni-directional application of equations; this
concept is borrowed from rewriting theory. The confluence condition is
satisfied, for instance, in the case where the monoid is generated by a
category whose morphisms act as partial maps on the space.

The basic construction of the globalization works in many topological
categories; here, we concentrate on topological spaces on the one hand,
and metric spaces on the other hand. For the topological case, we prove
that, under confluence, the original space is topologically embedded in
its globalization (and we provide an example which shows that this
result fails in the non-confluent case). Moreover, we show that the
globalization inherits normality and dimension from the original
space. Since free homogeneous extensions are globalizations for
(confluent) `singleton partial actions', this entails the corresponding
results for such extensions.

The metric setting is best considered in the larger category of
pseudometric spaces. Requiring confluence throughout, we prove an
embedding theorem, and we show that for an important class of cases, the
pseudometric globalization and the metric globalization coincide. We
demonstrate that, in these cases, dimension is preserved. Furthermore,
under suitable compactness assumptions, we prove existence of geodesic
paths; by consequence, the globalization of a path metric
space~\cite{Gromov99} is again a path metric space.

For every metric space, there exists an isometric embedding into a
metrically ultrahomogeneous space of the same weight. This is a part of
a recent result by Uspenskij~\cite{Uspenskij98}, and well-known for the
case of separable spaces~\cite{Urysohn27} (see also \cite{Uspenskij90};
for further information about Urysohn spaces, see
\cite{CGK,Gromov99,Pestov99,Vershik98}).  We show that in many cases the
metric globalization preserves the dimension. This implies that every
metric space $X$ admits a closed isometric embedding into an
ultrahomogeneous metric space $Z$ of the same dimension and
cardinality. It is an open question if
$Z$ can be chosen in such a way that the weight of $X$ is also
preserved.

\section{Confluently Generated Monoids}\label{free-monoid-section}

In preparation for the central notion of `well-behaved' partial action,
we now introduce a class of monoid presentations for which the word
problem is solvable by means of head-on application of directed
equations, i.e.\ by the classical rewriting method as used, up to now,
mainly in computer science applications such as $\lambda$-calculus and
automatic theorem proving~\cite{BaaderNipkow98,Klop92} (see however
\cite{SchroderHerrlich00a,SchroderHerrlich01b} for applications to
extensions of categories).

We recall that a \emph{monoid presentation} $\pres{G}{R}$ consists of a
set $G$ of \emph{generators} and a relation $R\subset G^*\times G^*$,
where $G^*$ is the set of \emph{words} over $G$, i.e.\
$G^*=\bigcup_{n=0}^\infty G^ n$. Here, we explicitly insist that $R$ is
a directed relation (rather than symmetric); the elements $(l,r)$ of
$R$, written $l\redstep r$, are called \emph{reduction rules} with
\emph{left side} $l$ and \emph{right side} $r$. Words are written either
in the form $(g_n,\dots,g_1)$ or, where this is unlikely to cause
confusion, simply in the form $g_n\dots g_1$. One way of describing the
monoid engendered by $\pres{G}{R}$ is as follows. The set $G^*$ is made
into a monoid by taking concatenation of words as multiplication,
denoted as usual simply by juxtaposition; the unit is the empty word
$()$. From $R$, we obtain a \emph{one-step reduction} relation
$\redstep$ on $G^*\times G^*$ by putting $w_1lw_2\redstep w_1rw_2$
whenever $(l,r)\in R$ and $w_1,w_2\in G^*$. Let $\redeq$ denote the
equivalence relation generated by $\redstep$; then the monoid $M$
described by $\pres{G}{R}$ is $G^*/\redeq$.

It is well known that the word problem for monoids, i.e.\ the question
whether or not $w_1\redeq w_2$ for given words $w_1$, $w_2$, is in
general undecidable. However, one can sometimes get a grip on the word
problem by means of normal forms: a word $w$ is called \emph{normal} if
it cannot be reduced under $\redstep$, i.e.\ if there is no word $w'$
such that $w\redstep w'$ (otherwise $w$ is called \emph{reducible});
thus, a word is normal iff it does not contain a left side of a
reduction rule. A normal word $w'$ is called a \emph{normal form} of a
word $w$ if $w\redeq w'$. We say that a monoid presentation is
\emph{noetherian} or \emph{well-founded} if the relation $\to$ is
well-founded, i.e.\ if there is no infinite sequence of reductions
$w_1\redstep w_2\redstep\dots$; this property guarantees existence, but
not uniqueness of normal forms. However, one can characterize those
cases where one does have uniqueness of normal forms. We denote the
transitive and reflexive closure of $\redstep$ by $\redge$ (reversely:
$\redle$); if $w\redge w'$, then $w'$ is said to be a \emph{reduct} of $w$.
\begin{myprop}\label{newman-prop} 
For a noetherian monoid presentation $\pres{G}{R}$, the following are
equivalent:
\begin{rmenumerate}
\item Each word in $G^*$ has a unique normal form.
\item Each word in $G^*$ has a unique normal reduct.
\item Whenever $w\redge s_1$ and $w\redge s_2$, then there exists a
  \emph{common reduct} $t\in G^*$ of $(s_1,s_2)$, i.e.\ $s_1\redge t$
  and $s_2\redge t$.
\item Whenever $w\redstep s_1$ and $w\redstep s_2$, then there is a
	common reduct of $(s_1,s_2)$.
\end{rmenumerate}
\end{myprop}
This proposition is a special case of a central lemma of rewriting
theory often referred to as \emph{Newman's Lemma} (see e.g.\
\cite{Klop92}, Theorem~1.0.7.). Condition (iii) is called
\emph{confluence}, while Condition (iv) is called \emph{weak
confluence}. The importance of the criterion lies in the fact that weak
confluence is often reasonably easy to verify. In particular, it is
enough to check weak confluence for so-called \emph{critical pairs},
i.e.\ cases where left sides of reductions rules overlap. More
precisely,
\begin{quote}
\em one can restrict Condition (iv) to words $w$ that are completely
made up of the overlapping left sides of the two involved reduction
rules
\end{quote}
(including the case that one of these left sides is contained in the
other); it is easy to see that this restricted condition is equivalent
to the original Condition (iv). Since the proof of
Proposition~\ref{newman-prop} is both short and instructive, we repeat
it here:
\begin{pf}
\emph{(i) $\implies$ (iv):} By the noetherian property, there exist
normal words $t_1$ and $t_2$ such that $s_1\redge t_1$ and $s_2\redge
t_2$. Then $t_1$ and $t_2$ are normal forms of $w$. By (i), we conclude
$t_1=t_2$.

\emph{(iv) $\implies$ (iii):} We proceed by the principle of
\emph{noetherian} or \emph{well-founded induction}, i.e.\ we prove the
claim for $w$ under the assumption that it holds for all proper reducts
of $w$. We can assume w.l.o.g.\ that both $w\redge s_1$ and $w\redge
s_2$ involve at least one reduction step, i.e.\ we have $w\redstep
w'_1\redge s_1$ and $w\redstep w'_2\redge s_2$. By (iv), we obtain a
common reduct $t$ of $(w'_1,w'_2)$. By the inductive assumption, we
obtain common reducts $r_1$ of $(s_1,t)$ and $r_2$ of $(s_2,t)$; again
by the inductive assumption, there is a common reduct of $(r_1,r_2)$,
which is then also a common reduct of $(s_1,s_2)$.

\emph{(iii) $\implies$ (ii):} Existence follows immediately from the
noetherian property. Concerning uniqueness, just observe that the
existence of a common reduct of two normal words implies their equality.

\emph{(ii) $\implies$ (i):} Whenever $w\redstep w'$, then (ii) implies
that $w$ and $w'$ have the same normal reduct. Thus, since $\redeq$ is
the equivalence relation generated by $\redstep$, this holds also
whenever $w\redeq w'$. In particular, for normal words $w$ and $w'$,
$w\redeq w'$ implies $w=w'$.
\rightqed
\end{pf}

\begin{rmdefn}\label{def:confluent-pres}
A noetherian monoid presentation is called \emph{confluent} if it
satisfies the equivalent conditions of Proposition~\ref{newman-prop} and
does not contain reduction rules with left side $g$, where $g\in G$.
\end{rmdefn}
The requirement that there are no left sides consisting of a single
generator can be satisfied for any noetherian monoid presentation by
removing superfluous generators, since for a reduction rule with left
side $g$, the noetherian condition implies that $g$ cannot occur on the
right side. Moreover, a noetherian monoid presentation cannot contain a
reduction rule with left side $()$. Thus, in confluent monoid
presentations any word with at most one letter is normal.

\begin{rmexmp}\label{exmp:confluentpres}
\begin{rmenumerate}
\item 	Every monoid has a trivial confluent presentation: take all 
	elements as generators, with reduction rules $uv\redstep p$ whenever
	$uv=p$. 
\item	The free monoid over a set $G$ of generators trivially has a
	confluent presentation $\pres{G}{\emptyset}$.
\item	The free group over a set $S$ of generators, seen as a monoid,
	has a confluent presentation $\pres{S\oplus S^{-1}}{R}$, where
	$\oplus$ denotes the disjoint union and
	$R$ consists of the reduction rules $ss^{-1}\redstep e$,
	$s^{-1}s\redstep e$ for each $s\in S$.
\item The free product $M_1*M_2$ of two monoids $M_1,M_2$ with confluent
	presentations $\pres{G_i}{R_i}$, $i=1,2$, respectively, has a
	confluent presentation $\pres{G_1\oplus G_2}{R_1\oplus R_2}$. If
	$M_1$ and $M_2$ are groups, then $M_1*M_2$ is a group, the
	free product of $M_1$ and $M_2$ as groups.
\item	The product $M_1\times M_2$ of two monoids $M_1$, $M_2$ with
	confluent presentations $\pres{G_i}{R_i}$, $i=1,2$,
	respectively, has a confluent presentation $\pres{G_1\oplus
	G_2}{R}$, where $R$ consists of all reduction rules in $R_1$ and $R_2$
	and the additional reduction rules $gh \redstep hg$ whenever $g\in
	G_2$, $h\in G_1$.
\item Given a subset $A$ of a monoid $M$ that consists of left
	cancellable elements, the monoid $M_A$ obtained by freely
	adjoining left inverses for the elements of $A$ has a confluent
	presentation $\pres{G}{R}$ as follows: we can assume that none
	of the elements of $A$ has a right inverse (since a right
	inverse of a left cancellable element is already a left
	inverse). Then  $G$ consists of the elements of $M$ and a new
	element $l_a$ for each $a\in A$; $R$ consists of the
	reduction rules for $M$ according to (i) and the reduction rules
	$(l_a,au)\redstep (u)$ for each $a\in A$, $u\in M$. This is a
	special case of a construction for categories discussed
	in~\cite{SchroderHerrlich00a}.
\item The infinite dihedral group has a confluent presentation
	$\pres{\{a,b,b^{-1}\}}{R}$, where $R$ consists of the reduction
	rules $bb^{-1}\redstep e$, $b^{-1}b\redstep e$, $aa\redstep e$,
	$ab\redstep b^{-1}a$, and $ab^{-1}\redstep ba$. (If the last
	reduction rule is left out, one still has a presentation of the
	same group, which however fails to be confluent.)
\item \label{exmp:catmonoid}
        Given a category $\BC$ \cite{AdamekHerrlich90,MacLane}, the
        monoid $M(\BC)$ induced by identifying all objects of $\BC$ (see
        e.g.~\cite{Borger77}) has a presentation $\pres{G}{R}$ given as
        follows. The set $G$ of generators consists of all morphisms of
        $\BC$. There are two types of reduction rules: on the one hand,
        rules of the form $(f,g)\redstep (f\circ g)$ for all pairs
        $(f,g)$ of composable morphisms in $\BC$, and on the other hand
        rules of the form $(\id_C)\redstep()$ for all objects $C$ of
        $\BC$. This presentation satisfies the conditions of
        Proposition~\ref{newman-prop}; it is turned into a confluent
        presentation in the stricter sense of
        Definition~\ref{def:confluent-pres} by removing all identities
        from the set of generators and modifying the reduction rule
        associated to a pair $(f,g)$ of morphisms to be
        $(f,g)\redstep()$ in case $f\circ g=\id$. This is a special case
        of the semicategory method introduced
        in~\cite{SchroderHerrlich00a}.
\end{rmenumerate}
\end{rmexmp}
Henceforth, we shall mostly denote elements of the monoid $M$ presented by
$\pres{G}{R}$ directly as words (or composites of letters) rather than
cluttering the notation by actually writing down equivalence classes of
words. E.g., phrases such as `$u$ has normal form $g_n\dots g_1$' means
that an element $u\in M$ is represented by the normal word
$(g_n,\dots,g_1)\in G^*$. The unit element will be denoted by $e$.

\begin{rmdefn}\label{def:prefix}
Let $M$ be a monoid with confluent presentation $\pres{G}{R}$.  An
element $u\in M$ with normal form $g_n\dots g_1$, where $g_i\in G$ for
$i=1,\dots,n$, is said to have \emph{length} $\length(u)=n$ (in
particular, $\length(e)=0$).  For a further $v\in M$ with normal form
$v=h_m\dots h_1$, we say that $uv$ is normal if $g_n\dots g_1h_m\dots
h_1$ is normal. We denote the order on $M$ induced by the prefix order
on normal forms by $\prefixle$; explicitly: we write $u \prefixle p$ iff
there exists $v$ such that $p=uv$ is normal. If additionally $u\neq p$,
then we write $u\prefixlt p$. The direct predecessor $g_n\dots g_2$ of
$u$ w.r.t.\ this order is denoted $\imprefix{u}$.
\end{rmdefn}

\section{Partial Actions and Globalizations}\label{preact-section}

Partial actions of groups have been defined and shown to coincide with
the restrictions of group actions to arbitrary subsets
in~\cite{Exel98}. We recall the definition, rephrased according
to~\cite{KellendonkLawson01}:
\begin{rmdefn}\label{def:grouppa}
Let $G$ be a group with unit $e$, let $X$ be a set, and let $\alpha$ be
a partial map $G\times X\to X$. We denote $\alpha(u,x)$ by $u\actson x$, 
with $\actson$ being right associative; i.e.\ $u\actson v\actson x$ denotes
$u\actson (v\actson x)$. The map $\alpha$ is called a \emph{partial
action} of $G$ on $X$ if, for each $x\in X$,
\begin{rmenumerate}
\item $e\actson x=x$,
\item if $u\actson x$ is defined for $u\in G$, then
  $u^{-1}\actson u\actson x=x$, and
\item if $u\actson v\actson x$ is defined, then 
  $(uv)\actson x=u\actson v\actson x$.
\end{rmenumerate}
Here, equality is to be read as \emph{strong} or \emph{Kleene} equality,
i.e.\ whenever one side is defined, then so is the other and the two
sides are equal. 
\end{rmdefn}
Concrete examples of partial group actions, including partial actions of
groups of M\"obius transforms, as well as further references can be
found in~\cite{KellendonkLawson01}. 
\begin{rmrem}
In~\cite{KellendonkLawson01}, partial actions are defined by
Conditions~(ii) and~(iii) above, and partial actions satisfying
Condition~(i) are called~\emph{unital}. The original definition of
partial actions~\cite{Exel98} includes Condition~(i).
\end{rmrem}
We generalize this definition to monoids as follows.
\begin{rmdefn}\label{def:monoidpa}
Given a set $X$, a \emph{partial action} of a monoid $M$ with unit $e$
on $X$ is a partial map
\begin{displaymath}
\alpha: M\times X\to X,
\end{displaymath}
with the notation $\alpha(u,x)=u\actson x$ as in Definition~\ref{def:grouppa}, such
that
\begin{rmenumerate}
\item	$e\actson x=x$ for all $x$, and
\item 	$(uv)\actson x = u\actson v\actson x$ whenever $v\actson x$  
	is defined.
\end{rmenumerate}
(Again, (ii) is a strong equation.)  Given two such partial actions of
$M$ on sets $X_1$, $X_2$, a map $f:X_1\to X_2$ is called
\emph{equivariant} if $u\actson f(x)$ is defined and equal to
$f(u\actson x)$ whenever $u\actson x$ is defined.
\end{rmdefn}
We explicitly record the fact that partial monoid actions indeed
generalize partial group actions:
\begin{myprop}
The partial monoid actions of a group $G$ are precisely its partial
group actions.
\end{myprop}
\begin{pf}
In the notation as above, let $e\actson x=x$ for all $x\in X$. We have
to show that Conditions~(ii) and~(iii) of Definition~\ref{def:grouppa}
hold iff Condition~(ii) of Definition~\ref{def:monoidpa} holds.

\noindent\emph{`If':} Condition~(iii) is immediate, since definedness of
$u\actson v\actson x$ entails definedness of $v\actson
x$. Moreover, if $u\actson x$ is defined, then by
Definition~\ref{def:monoidpa}~(ii), we have $u^{-1}\actson u\actson
x=(u^{-1}u)\actson x=e\actson x=x$; this establishes
Definition~\ref{def:grouppa}~(ii). 

\noindent\emph{`Only If':} The right-to-left direction of the strong
equation in Definition~\ref{def:monoidpa}~(ii) is just
Definition~\ref{def:grouppa}~(iii). To see the converse direction, let
$u,v\in G$, and let $v\actson x$ and $(uv)\actson x$ be defined; we have
to show that $u\actson v\actson x$ is defined. By
Definition~\ref{def:grouppa}~(ii), $v^{-1}\actson v\actson x=x$, so that
$(uv)\actson v^{-1}\actson v\actson x$ is defined; by
Definition~\ref{def:grouppa}~(iii), it follows that $(uvv^{-1})\actson
v\actson x$ is defined, and this is $u\actson v\actson x$.
\rightqed
\end{pf}

A partial action is equivalently determined by the partial maps 
\begin{displaymath}
\begin{array}{rccc}
u: & X & \to & X\\
   & x & \mapsto & u\actson x
\end{array}
\end{displaymath}
associated to $u\in M$.  The domain of $u:X\to X$ is denoted
$\domain(u)$.

Here, we are interested mainly in partial actions on spaces of some
kind. E.g., we call a partial action of $M$ on a topological space $X$
\emph{continuous} if the associated partial map $\alpha: M\times X\to X$
is continuous on its domain, where $M$ carries the discrete topology,
equivalently: if each of the maps $u:X\to X$ is continuous on
$\domain(u)$. A partial action is called \emph{closed} (\emph{open}) if
$\domain(u)$ is closed (open) for each $u\in M$, and \emph{strongly
closed} (\emph{strongly open}) if, moreover, $u:X\to X$ is
closed (open) on $\domain(u)$ for each $u$.

It is clear that a (total) action of $M$ on a set $Y$ induces a partial
action on each subset $X\subset Y$. This statement has a converse:

\begin{rmdefn}
Given a partial action of $M$ on $X$, its \emph{(universal)
globalization} consists of a set $Y$ with a total action of $M$ and an
equivariant map $i:X\to Y$ such that every equivariant map from $X$ to a
total action of $M$ factors uniquely through $i$. 
\end{rmdefn}
(Topological and metric globalizations are defined analogously,
requiring continuity and non-expansiveness, respectively, for all
involved maps.)

The globalization is easy to construct at the set level: the set $Y$ is
the quotient of $M\times X$ modulo the equivalence relation $\simeq$
generated by
\begin{equation}\label{eqn:generateeq}
(uv,x)\sim(u,v\actson x)\quad\textrm{whenever $v\actson x$ is defined}
\end{equation}
(the generating relation $\sim$ is reflexive and transitive, but unlike
in the case of groups fails to be symmetric). We denote the equivalence
class of $(u,x)$ by $[u,x]$. The action of $M$ is defined by $u\actson
[v,x]=[uv,x]$. Moreover, $i(x)=[e,x]$. This map makes $X$ a subset of
$Y$:

\begin{myprop}
The map $i:X\to Y$ defined above is injective, and the action of $M$ on
$Y$ induces the original partial action on $X$.
\end{myprop}
\begin{pf}
Define an equivalence relation $\rho$ on $M\times X$ by 
\begin{displaymath}
(u,x)\rho (v,y)\iff u\actson x=v\actson y,
\end{displaymath}
where again equality is strong equality. By
Definition~\ref{def:monoidpa}~(ii), $\rho$ contains the relation $\sim$
defined in Formula~(\ref{eqn:generateeq}) above. Thus, $\rho$ contains
also the equivalence $\simeq$ generated by $\sim$; i.e.\
$(u,x)\simeq(v,y)$ implies the strong equation $u\actson x=v\actson
y$. In particular, $(e,x)\simeq(e,y)$ implies $x=e\actson x=e\actson
y=y$, so that $i$ is injective. Moreover, it follows that
$(u,x)\simeq(e,y)$ implies that $u\actson x=y$ is defined, i.e.\ the
restriction of the action on $Y$ to $X$ is the given partial
action.\rightqed
\end{pf}
Thus, partial actions of monoids are precisely the restrictions of total
actions to arbitrary subsets. From now on, we will identify $X$ with
$i(X)$ whenever convenient. By the second part of the above proposition,
overloading the notation $u\actson x$ to denote both the action on $Y$
and the partial action on $X$ is unlikely to cause any confusion.

The proof of the above proposition shows that equivalence classes of
elements of $X$ are easy to describe; however, a similarly convenient
description is not generally available for equivalence classes of
arbitrary $(u,x)$ --- that is, $(u,x)\simeq (v,y)$ may mean that one has
to take a `zig-zag path' from $(u,x)$ to $(v,y)$ that uses the
generating relation $\sim$ of Formula~(\ref{eqn:generateeq}) both from
left to right and from right to left. However, the situation is better
for partial actions that have well-behaved presentations in the same
spirit as confluently presented monoids.

Let $\alpha$ be a partial action of a monoid $M$ on $X$, and let
$\pres{G}{R}$ be a confluent presentation of $M$. Then we regard the
restriction of $\alpha$ to $G\times X$ as a collection of additional
\emph{reduction rules}, i.e.\ we write
\begin{equation}\label{eqn:gen-reduct}
(g,x) \redstep (g\actson x)\quad\textrm{whenever $g\actson x$ is defined
for $g\in G$, $x\in X$},
\end{equation}
in addition to the reduction rules already given by $R$.  In the same
way as for monoid presentations, this gives rise to a \emph{one-step
reduction} relation $\redstep$ on the set $G^*\times X$, whose elements
we denote in either of the two forms $(g_n,\dots,g_1,x)$ or $g_n\dots
g_1\actson x$. Explicitly, we write
$(g_n,\dots,g_2,g_1,x)\redstep(g_n,\dots,g_2,g_1\actson x)$ whenever
$g_1\actson x$ is defined, and $w_1\actson x\redstep w_2\actson x$
whenever $w_1\redstep w_2$ for words $w_1,w_2\in G^*$. Moreover, we
denote the transitive and reflexive hull of $\redstep$ and the
equivalence relation generated by $\redstep$ on $G^*\times X$ by
$\redge$ and $\redeq$, respectively, and we use the terms \emph{normal},
\emph{normal form}, \emph{reduct}, and \emph{common reduct} as
introduced for words in $G^*$ in the previous section with the obvious
analogous meanings for words in $G^*\times X$. Since the additional
reduction rules always reduce the word length by $1$, it is clear that
reduction in $G^*\times X$ is also \emph{well-founded} (or
\emph{noetherian}), i.e.\ that there are no infinite reduction sequences
in $G^*\times X$. Thus, we have an analogue of
Proposition~\ref{newman-prop} (with almost literally the same proof):
\begin{myprop}\label{newman2-prop} 
In the above notation, the following are equivalent:
\begin{rmenumerate}
\item Each word in $G^*\times X$ has a unique normal form.
\item Each word in $G^*\times X$ has a unique normal reduct.
\item Whenever $w\redge s_1$ and $w\redge s_1$ in $G^*\times X$, then
  there exists a common reduct of $(s_1,s_2)$,
\item Whenever $w\redstep s_1$ and $w\redstep s_2$ in $G^*\times X$,
	then there exists a common reduct of $(s_1,s_2)$.
\end{rmenumerate}
\end{myprop}
In fact, the point behind all these analogies is that $(G^*\times
X,\redstep)$ is just another example of a rewrite system, and the above
proposition is another special case of Newman's Lemma. Concerning the
verification of \emph{weak confluence}, i.e.\ Condition (iv) above, we
remark that, besides checking confluence of $\pres{G}{R}$, it suffices
to consider cases of the form $w=g_n\dots g_1\actson x$, where $g_n\dots
g_1$ is the left side of a reduction rule in $R$ and $g_1\actson x$ is
defined.

\begin{rmdefn}\label{def:confluentPA}
A partial action of a monoid $M$ on a set $X$ is called \emph{confluent}
if $M$ has a confluent presentation $\pres{G}{R}$ (cf.\
Section~\ref{free-monoid-section}) such that the equivalent conditions
of Proposition~\ref{newman2-prop} hold for the associated reduction
relation $\redstep$ on $G^*\times X$, and such that this reduction
relation \emph{generates} the given partial action. The latter means
explicitly that, for $g_n\dots g_1\in G^*$,
\begin{displaymath}
(g_n\dots g_1)\actson x=y\quad\textrm{implies}\quad(g_n,\dots,g_1,x)\redge (y)
\end{displaymath}
(the converse implication holds by the definition of partial actions).

\emph{For the sake of brevity, we shall fix the notation introduced so
far ($\alpha$ for the action, $X$ for the space, $Y$ for the
globalization, $G$ for the set of generators etc.) throughout.}
\end{rmdefn}

By the generation condition, the quotient of $G^*\times X$ modulo the
equivalence relation $\redeq$ is the universal globalization constructed
above, so that we now have a way of deciding equivalence of
representations for elements of the globalization outside $X$, namely
via reduction to normal form.  This will allow us to reach a good
understanding of the properties of the globalization as a space.

In typical applications, a confluent partial action will often be given
in terms of a monoid presentation $\pres{G}{R}$ and a partial map
$G\times X\to X$; in this case, the partial action of the monoid $M$
presented by $\pres{G}{R}$ is \emph{defined} by putting $g_n\dots
g_1\actson x=y\iff (g_n,\dots,g_1,x)\redge (y)$. Verifying the
conditions of Proposition~\ref{newman2-prop} then guarantees that this
does indeed define a partial action.

\begin{rmexmp}\label{exmp:confluentpa}
\begin{rmenumerate}
\item A partial action of $M$ is confluent w.r.t.\ the trivial confluent
	presentation of $M$ (cf.\ Example~\ref{exmp:confluentpres}~(i))
	iff, whenever $v\actson x$, then either $(uv)\actson x$ is
	defined or $(u,v\actson x)= (uv,x)$: to see this, assume
	$(u,v\actson x)\neq (uv,x)$; then $(uv)\actson x$ is the only
	possible common reduct of the reducts $(u,v\actson x)$ and
	$(uv,x)$ of $(u,v,x)$. Most of the time, this is a rather too
	strong property to require. In particular, if $M$ is a group,
	then this holds iff, for each $v\neq e$, definedness of
	$v\actson x$ implies definedness of $(uv)\actson x$ for each $u$
	--- this means that the partial action at hand is essentially
	just a total action on the subset $\{x\mid v\actson x\textrm{ is
	defined for some }v\neq e\}$ of $X$.
\item	Partial actions of the free monoid over $G$ are always confluent
	w.r.t.\ the confluent presentation $\pres{G}{\emptyset}$.
\item \label{exmp:freegrouppa} Partial actions of the free group over $S$
	are always confluent w.r.t.\ the confluent presentation of
	Example~\ref{exmp:confluentpres}~(iii).
\item	Two confluent partial actions of monoids $M_1$ and $M_2$ on a
	set $X$, respectively, give rise to a confluent partial action of 
	$M_1*M_2$ on $X$ w.r.t.\ the confluent presentation
	given in Example~\ref{exmp:confluentpres}~(iv).
\item	A total action of $M$ on $X$ can be extended to a confluent partial
	action on $X$ of the extended monoid $M_A$ of
	Example~\ref{exmp:confluentpres}~(vi) w.r.t.\ the confluent
	presentation given there (by putting
	$l_a\actson(au\actson x)=u\actson x$ for each $a\in A$, $u\in
	M$, $x\in X$) iff each $a\in A$ acts injectively on $X$.
\item	A partial action of the infinite dihedral group is confluent
	w.r.t.\ the confluent presentation given in
	Example~\ref{exmp:confluentpres}~(vii) iff
	\begin{alenumerate}
	\item	$a\actson x$ and $ab\actson x$ are defined whenever
		$b\actson x$ is defined, and
	\item	$a\actson x$ and $ab^{-1}\actson x$ are defined whenever
		$b^{-1}\actson x$ is defined.
	\end{alenumerate}
\item \label{exmp:preact} A partial action of the monoid $M(\BC)$
  generated by a small category $\BC$ as in
  Example~\ref{exmp:confluentpres}~\romanref{exmp:catmonoid} on a set
  $X$ is confluent (w.r.t.\ the given confluent presentation of
  $M(\BC)$) iff, whenever $f$ and $g$ are composable morphisms in $\BC$
  and $g\actson x$ is defined, then either $(f\circ g,x)=(f,g\actson x)$,
  or $(f\circ g)\actson x$ is defined (and hence also $f\actson
  (g\actson x)$). 

  In particular, this is the case if the partial action
  is given by a functor from $\BC$ into the category $\BS(X)$ of maps
  between subsets of $X$; this generalizes the preactions of groupoids
  considered in~ \cite{Megrelishvili85,Megrelishvili86,Megrelishvili00}.
  Here, we need only the simpler case that $\BC$ is actually a
  subcategory of $\BS(X)$. Explicitly, such a subcategory determines a
  confluent partial action of $M(\BC)$ as follows: if $f:A\to B$ is a
  morphism of $\BC$, i.e.\ a map between subsets $A$ and $B$ of $X$,
  then $f\actson x$ is defined iff $x\in A$, and in this case equal to
  $f(x)$. Analogously, one obtains a continuous partial action on a
  topological space $X$ from a subcategory of the category $\BT(X)$ of
  continous maps between subspaces of $X$ etc.
\end{rmenumerate}
\end{rmexmp}

\begin{rmrem}
Due to Example~\ref{exmp:confluentpres}~(i), it does not make sense to
regard the existence of a confluent presentation as a property of a
monoid; rather, a confluent presentation is considered as extra
structure on a monoid. Contrastingly, the results about confluent
partial actions presented below depend only on the existence of a
confluent presentation; in the few places where we do make reference to
the generating system in definitions, these definitions will turn out to
be in fact independent of the chosen generating system by virtue of
subsequently established results (see for example
Definition~\ref{def:degenerate} and
Proposition~\ref{prop:pseudometric}). Thus, we mostly think of confluence
of a partial action as a property; Example~\ref{exp:non-confluent} will
show that not all partial actions have this property.
\end{rmrem}

As in the case of monoids, we usually denote the elements of $Y$
directly by their representatives in $G^*\times X$ rather than as
explicit equivalence classes. Of course, we can still represent elements
of $Y$ as pairs $(u,x)\in M\times X$. We will say that $(u,x)$ or
$u\actson x$ is in normal form if $g_n\dots g_1\actson x$ is in normal
form, where $g_n\dots g_1$ is the normal form of $u$; similarly, we
write $u\actson x\redge v\actson y$ if this relation holds with $u$ and
$v$ replaced by their normal forms, etc. By the definition of confluent
presentation, $g\actson x$ is normal for $g\in G$, $x\in X$, whenever
$g\actson x$ is undefined in $X$. Moreover, $e\actson x$ is always
normal. We put
\begin{displaymath}
R_u =\{x \in X\mid   u\actson x\textrm{ is normal} \}=X\setminus \domain(g_1),
\end{displaymath}
where $u$ has normal form $g_n\dots g_1$ (note that $R_e=X$). The action
of $u$ gives rise to a bijective map $u: R_u \to u\actson R_u$.

\begin{rmdefn}
An element $\globa\in Y$ with normal form $g_n\dots g_1\actson
x$ is said to have \emph{length} $\length(\globa)=n$. We put
\begin{displaymath}
Y_n =\{\globa \in Y\mid   \length(\globa) \leq n \}.
\end{displaymath}
\end{rmdefn}

Of course, a confluent partial action is continuous iff the partial map
$g:X\to X$ is continuous for each generator $g\in G$. A similar
reduction holds for the domain conditions (closedness etc.); cf.\
Section~\ref{sec:globalizations}.

We finish this section by exhibiting an example of a partial action that
fails to be confluent. This relies on an observation
concerning the structure of the universal globalization $Y$ of a
confluent partial action.

\eat{
\begin{rmdefn}\label{def:cycle}
A sequence $u_1,\dots,u_n$ of elements $u_i\in M$, indexed modulo $n$,
is called an \emph{$n$-cycle} of $\alpha$ if
\begin{displaymath}
u_i\actson X\cap u_{i+1}\actson X\neq\emptyset
\end{displaymath}
in $Y$ for each $i$. Such an $n$-cycle is called \emph{trivial} if there
exists $v\in M$ such that
\begin{displaymath}
u_i\actson X\cap u_{i+1}\actson X\subset v\actson X
\end{displaymath}
for each $i$.
\end{rmdefn}
}

\begin{mylem}\label{lem:red-pass}
Let $\alpha$ be confluent, and let $a=u\actson x$ have normal form
$v\actson y$. Then $a\in w\actson X$ whenever $v\prefixle w\prefixle u$
in the prefix order (cf.\ Definition~\ref{def:prefix}).
\end{mylem}
\begin{pf}
The reduction from $(u, x)$ to $(v, y)$ works by taking the normal form
of $u$ and then shifting letters from left to right according to 
Formula~(\ref{eqn:gen-reduct}). Thus, there must be an intermediate step
of the form $(w,z)$, which proves the claim.\rightqed
\end{pf}
\begin{myprop}\label{prop:triples}
If $\alpha$ is confluent, then for every triple $(u_1,u_2,u_3)\in M^3$
(indexed modulo $3$), there exists $w\in M$ such that, for $i=1,2,3$,
\begin{displaymath}
u_i\actson X\cap u_{i+1}\actson X\subset w\actson X
\end{displaymath}
in $Y$.
\end{myprop}
\begin{pf}
Let $w_i=u_i\meet u_{i+1}$ for $i=1,2,3$. Here, $\meet$ denotes the meet
in the prefix order (cf.\ Definition~\ref{def:prefix}), i.e. the
largest common prefix. Now since for each $i$, $w_i$ and $w_{i+1}$ are
both prefixes of $u_{i+1}$, they are comparable under the prefix order;
i.e.\ the $w_i$ form a chain. We can assume w.l.o.g.\ that $w_1$ is
the largest element of this chain.

Then $w:=w_1$ has the claimed property. Indeed, if $a=u_i\actson
x=u_{i+1}\actson y$, then by confluence, $a$ must have normal form
$a=v\actson z$, where $v\prefixle u_i$ and $v\prefixle u_{i+1}$. Thus,
$v\prefixle w_i\prefixle w_1$; by Lemma~\ref{lem:red-pass}, this implies
$a\in w_1\actson X$, because we have $w_1\prefixle u_i$ or $w_1\prefixle
u_{i+1}$.\rightqed
\end{pf}

\begin{rmexmp}\label{exp:non-confluent}
Let $\Klein$ denote the Klein four-group $\{e,u,v,uv\}$, and let
$\alpha$ be the partial action of $\Klein$ on the set $\{0,1,2\}$
defined by letting $u$, $v$, and $uv$ act as partial identities defined
on the domains $\{0\}$, $\{1\}$, and $\{2\}$, respectively. Then the
triple $(e,u,uv)\in\Klein^3$ violates the property in
Proposition~\ref{prop:triples}. To see this, we show that
\begin{displaymath}
0\notin v\actson X\cup uv\actson X,\quad 
u\actson 1\notin X,\quad\textrm{and}\quad
2\notin u\actson X. 
\end{displaymath}
The equivalence class of $(e,0)$ in $M\times X$ is $\{(e,0),(u,0)\}$,
because this set is closed under the generating relation $\sim$ of
Formula~(\ref{eqn:generateeq}) above, so that indeed $0\notin v\actson
X\cup uv\actson X$; the other claims are proved similarly. 
Since we have
\begin{displaymath}
0\in X\cap u\actson X,\quad u\actson 1\in u\actson X\cap uv\actson X,
\quad\textrm{and}\quad 2\in uv\actson X\cap X,
\end{displaymath}
we have shown that there is no  $w\in\Klein$ such that
$w\actson X$ contains all three pairwise intersections of $X$, $u\actson
X$, and $uv\actson X$. Thus, $\alpha$ fails to be confluent.
\end{rmexmp}

\section{Topological Globalizations}\label{sec:globalizations}

We now move on to discuss universal globalizations of continuous partial
actions of a monoid $M$ on a topological space $X$; here, the
universality is, of course, to be understood w.r.t.\ continuous
equivariant maps. The main result of this section states essentially
that globalizations of \emph{confluent} partial actions of monoids are
topological embeddings. A corresponding result for \emph{open} partial
group actions (without confluence) is established
in~\cite{KellendonkLawson01} and in~\cite{Abadie03}. We shall provide an
example that shows that the result fails for arbitrary partial group
actions.  \eat{As in several other globalization constructions
(see~\cite[ch. 9.8]{Beardon83}, \cite{Megrelishvili86,Megrelishvili00},
\cite[section 4]{Uspenskij98}, and references to further recent work in
\cite{KellendonkLawson01}), the global set of~\cite{KellendonkLawson01}
can be obtained as a quotient set of $G \times X$ under an equivalence
relation defined by a certain inverse semigroup of partial bijections.}

The universal globalization of a continuous partial action is
constructed by endowing the globalization $Y$ constructed above with the
final topology w.r.t.\ the maps
\begin{displaymath}
\begin{array}{rccc}
u: & X & \to & Y\\
   & x & \mapsto & u\actson x,
\end{array}
\end{displaymath}
where $u$ ranges over $M$ (i.e.\ $V\subset Y$ is open iff $u^{-1}[V]$ is
open in $X$ for each $u\in M$); equivalently, the topology on $Y$ is the
quotient topology induced by the map $M \times X\to Y$, where $M$
carries the discrete topology.  This ensures the desired universal
property: given a continuous equivariant map $f:X\to Z$, where $M$ acts
globally (and continuously) on $Z$, the desired factorization $f^\#:Y\to
Z$ exists uniquely as an equivariant map by the universal property of
$Y$ at the level of sets. In order to establish that $f^\#$ is
continuous, it suffices to show that $f^\#u:X\to Z$ is continuous for
each $u\in M$; but $f^\#u$ is, by equivariance of $f^\#$, the map
$x\mapsto u\actson f(x)$, hence continuous.

\eat{given a morphism $(G,g)$ (as
defined in Section~\ref{preact-section}) from the partial action $F$ to an
action $H$ of a monoid $N$ on a further space $Z$, the universal
property of $Y$ at the set level induces a morphism of monoid actions
$(G^\#,g^\#)$ from $Y$ to $Z$. Composed with $u:X\to Y$, where $u\in
M(\BC)$, $g^\#$ becomes the continuous map $x\mapsto G(u)\actson g(x)$
(since $g^\#(u\actson x)=G(u)\actson g^\#(x)=G(u)\actson g(x)$); hence,
$g^\#$ is continuous.}

Under additional assumptions concerning the domains, the inclusion
$X\into Y$ is extremely well-behaved:
\begin{myprop} If $\alpha$ is closed (open),
then the map $X\into Y$ is closed (open), in particular a topological
embedding. 
\end{myprop}
(The open case for partial group actions appears in
\cite{Abadie03,KellendonkLawson01}.)
\begin{pf}
Let $A\subset X$ be closed (open). Then $u^{-1}[A]$ is closed (open) in
$\domain(u)$ and hence in $X$ for each $u\in M$; thus, $A$ is closed
(open) in $Y$.
\rightqed
\end{pf}
The embedding property fails in the general case:
\begin{rmexmp}\label{exmp:non-embedding} 
We proceed similarly as in Example~\ref{exp:non-confluent}. Let $\Klein$
denote the Klein four-group $\{e,u,v,uv\}$, and let $\alpha$ be the
partial action of $\Klein$ on the closed interval $X=[-1,1]$ defined as
follows: let $\domain(u)=A=\left\{\frac{1}{2}\right\}$, let
$\domain(v)=B=\left\{\frac{1}{n}+\frac{1}{2}\mid n\in\Nat, n\ge
2\right\}$, and let $\domain uv=C=[-1,1]\cap\Rat$. Let $u$ and $v$ act
as the identity on $A$ and $B$, respectively, and let $(uv)\actson x=
-x$ for $x\in C$. It is easily checked that $\alpha$ is indeed a partial
group action. As in Example~\ref{exp:non-confluent}, one shows that
$\alpha$ fails to be confluent, because the triple $(1,u,uv)$ violates
the property in Proposition~\ref{prop:triples} (alternatively,
non-confluence of $\alpha$ can be deduced from the following and
Corollary~\ref{cor:top-embed}).

We claim that \emph{the globalization $X\into Y$ of $\alpha$ fails to
be a topological embedding} (which, incidentally, implies that $Y$ fails
to be Hausdorff, since $X$ is compact and $X\into Y$ is injective). To
see this, let $U$ be the open set $(0,1)$ in $X$. We show that $U$ fails
to be open in $Y$, i.e.\ that $V\cap X\neq U$ for each open $V\subset Y$
such that $U\subset V$; in fact, such a $V$ always contains a negative
number:

We have $u\actson \frac{1}{2}=\frac{1}{2}\in V$, i.e.\ $\frac{1}{2}\in
u^{-1}[V]$. Therefore the open set $u^{-1}[V]\subset X$ intersects $B$,
i.e.\ we have $b\in B$ such that $(uv)\actson b=u\actson b\in V$. Thus,
the open set $(uv)^{-1}[V]\subset X$ intersects $C\cap(0,1]$, so that we
obtain $c\in C\cap(0,1]$ such that $(uv)\actson c\in V$; but then
$(uv)\actson c=-c$ is a negative number.
\end{rmexmp}
Notice that it is not possible to repair the embedding property by just
changing the topology on $Y$: the topology is already as large as
possible (being a final lift of maps that are certainly expected to be
continuous), and the failure of $X\into Y$ to be an embedding is due to
$Y$ having \emph{too few} open sets. This pathology does not happen in
the confluent case:
\begin{mythm}\label{thm:top-embed}
If $\alpha$ is confluent, then the map $u:R_u\to Y$ (cf.\
Section~\ref{preact-section}) is a topological embedding for each $u\in
M$.
\end{mythm}
\begin{mycor}\label{cor:top-embed}
If $\alpha$ is confluent, then the globalization $X\into Y$ is a
topological embedding.
\end{mycor}
(It is unlikely that the converse holds, i.e.\ that that confluence is
also a necessary condition for $X\into Y$ to be an
embedding.)
\begin{mypf}{Corollary~\ref{cor:top-embed}} The inclusion
$X\into Y$ is the map $e:R_e\to Y$.
\rightqed
\end{mypf}

\begin{mypf}{Theorem~\ref{thm:top-embed}}
All that remains to be shown is that the original topology of $R_u$
agrees with the subspace topology on $u\actson R_u$ w.r.t.\ $Y$, i.e.\
that, whenever $U$ is open in $R_u$, then there exists an open $\bar
U\subset Y$ such that $\bar U\cap u\actson R_u=u\actson U$.

We define $\bar U$ as the union of a system of subsets $U_v\subset Y$ to
be constructed below, indexed over all $v\in M$ such that $u\prefixle v$
(this is the prefix ordering of Definition~\ref{def:prefix}, which
depends on confluence. As announced above, we reuse notation without
further comments.), with the following properties for each $v\prefixge u$:
\begin{rmenumerate}
\item	$U_p\subset U_v$ whenever $u\prefixle p\prefixle v$.
\item	$U_v\cap u\actson R_u=u\actson U$.
\item	$v^{-1}[U_v]$ is open in $X$.
\item	Each $\globa\in U_v \setminus U_{\imprefix{v}}$ has normal form 
	$v\actson x$ for some $x$.
\end{rmenumerate}
Then certainly 
\begin{displaymath}
\bar U\cap u\actson R_u=u\actson U.
\end{displaymath}
Moreover, the properties above imply
\begin{rmenumerate}
\setcounter{enumi}{4}
\item   For each $v\in M$, $v\actson x\in\bar U$ implies $u\prefixle v$
        and $v\actson x\in U_v$.
\end{rmenumerate}
To prove (v), let $p$ be the minimal $p\prefixge u$ w.r.t.\ $\prefixle$
such that $v\actson x\in U_p$. By (iv), $v\actson x$ has normal form
$p\actson y$ for some $y$, so that $p\prefixle v$, and hence in
particular $u\prefixle v$. By (i), we obtain $v\actson x\in U_v$ as
required. Now (v) enables us to show that $\bar U$ is open: we have to
show that $v^{-1}[\bar U]$ is open for each $v\in M$. By (v), this set
is empty in case $u\not\prefixle v$. Otherwise, we have, again by (v),
\begin{displaymath}
v^{-1}[\bar U]=v^{-1}[U_v]
\end{displaymath}
which is open in $X$ by (iii).

The system $(U_v)$ is constructed by induction over the prefix order,
starting from $U_u=U$ (where `$U_{\imprefix{u}}$' is to be replaced by
$\emptyset$ in (iv)). Now let $v\in M$, where $u\prefixlt v$, have
normal form $v=g_n\dots g_1=\imprefix{v}g_1$, and assume that the $U_p$
are already constructed as required for $u\prefixle p\prefixlt v$. The
set
\begin{displaymath}
B=(\imprefix{v})^{-1}[U_{\imprefix{v}}]
\end{displaymath}
is open in $X$ by the inductive assumption. Thus, $g_1^{-1}[B]$ is open
in the domain $D\subset X$ of $g_1$, i.e. equal to $D\cap V$, where $V$
is open in $X$.  Let
\begin{displaymath}
C=V\setminus D.
\end{displaymath}
Note that, for $x\in C$, $v\actson x$ is normal. Now $U_v$ is defined
as
\begin{displaymath}
U_v=U_{\imprefix{v}}\cup v\actson C.
\end{displaymath}
It is clear that this definition satisfies (i), (ii), and (iv) above. In
order to verify (iii), let $x\in X$. Then $v\actson x$ is normal and in
$U_v$ iff $x\in C$. If $v\actson x$ is reducible, i.e.\ if $g_1(x)$ is
defined in $X$, then $v\actson x\in U_v$ iff $\imprefix{v}\actson
(g_1(x))\in U_{\imprefix{v}}$ iff $g_1(x)\in B$. Thus,
\begin{displaymath}
v^{-1}[U_v]=C\cup g_1^{-1}[B]=(V\setminus D)\cup (V\cap D)=V,
\end{displaymath}
which is open in $X$.
\rightqed
\end{mypf}

\begin{rmexmp}\label{exmp:shimrat}
A very basic example of a partial action on $X$ produces the free
homogeneous space over $X$, as follows. The full subcategory $\BC$ of
$\BT(X)$ spanned by the singleton subspaces induces a partial action as
described in Example~\ref{exmp:confluentpa}~\romanref{exmp:preact}. The
presentation of the monoid $M(\BC)$ generated by $\BC$ can be described
as follows: the generators are of the form $(xy)$, where $x,y\in X$ with
$x\neq y$, and the relations are $(xy)(yz)\redstep(xz)$ when $x\neq z$,
and $(xy)(yz)\redstep()$ otherwise (thus, one may leave out the brackets
and just write $xx=e$). The corresponding globalization is easily seen
to be homogeneous. There are known ways to produce this homogeneous
space, in particular \emph{Shimrat's construction}~\cite{Shimrat54} and
the construction given by Belnov~\cite{Belnov78}, who also establishes a
kind of universal property for the extension. It can be checked that the
spaces resulting from these constructions coincide with our
globalization in this special case (see~\cite{Megrelishvili85} for more
details).
\end{rmexmp}

\section{Preservation of Topological Properties}

We will now investigate how topological properties of a space are or are
not handed on to its globalization with respect to a continuous partial
action $\alpha$.
\begin{mythm}
\label{T1} 
If $\alpha$ is confluent and $X$ is a $T_1$-space, then $Y$ is $T_1$ iff
$u^{-1}[\{x\}]$ is closed in $X$ for each $u\in M$ and each $x\in X$.
\end{mythm}
\begin{pf}
The `only if' direction is immediate. In order to prove the `if'
direction, we have to show that the latter condition implies that
$u^{-1}[\{a\}]$ is closed in $X$ for each $a\in Y$. Let $a$ have normal
form $v\actson x$. Then $u\actson y=v\actson x$ for $y\in X$ iff we have
$u=vp$ normal and $p\actson y=x$, where $p$ is necessarily uniquely
determined. Thus, $u^{-1}[\{a\}]$ is the closed set $p^{-1}[\{x\}]$ in
$X$ if $v\prefixle u$; otherwise, $u^{-1}[\{a\}]$ is empty.
\rightqed
\end{pf}

There are many typical cases in which this necessary and sufficient
condition is easily seen to be satisfied, such as the following.
\begin{mycor}\label{cor:T1-closed}
If $X$ is $T_1$ and $\alpha$ is closed, then $Y$ is $T_1$.
\end{mycor}
\begin{mycor}\label{cor:T1-group}
If $X$ is $T_1$ and $M$ is a group, then $Y$ is $T_1$.
\end{mycor}
\begin{mycor}\label{cor:T1-finite-fibres}
If $X$ is $T_1$ and for each generator $g\in G$, the partial map $g:X\to
X$ has finite fibres, then $Y$ is $T_1$.
\end{mycor}
(The latter corollary includes the case that all generators act
injectively.)
\begin{mypf}{Corollary~\ref{cor:T1-finite-fibres}}
By induction over the length of $u\in M$, one shows that $u^{-1}[\{x\}]$
is finite and hence closed for each $x\in X$.
\rightqed
\end{mypf}

\eat{
If $G\subset M$ generates $M$, then elements $w=(g_n,\dots,g_1)$ of
$G^*$ act on $X$ as partial continuous maps $w^*$ in the obvious way:
for $x\in X$, let
\begin{displaymath}
w^*(x)=g_n \actson \dots \actson g_1 \actson x
\end{displaymath}
(where the expression on the right is understood to be undefined if the
application of one of the $g_i$ fails). Note that words $w,v$ that
represent the same element of $M$ need not satisfy $w^*=v^*$.
}

For confluent actions, the domain conditions introduced in
Section~\ref{preact-section} can be reduced to the generating set $G$:
\begin{myprop}\label{prop:closed}
Let $\alpha$ be confluent. Then $\alpha$ is closed (open) iff
$\domain(g)$ is closed (open) for each $g\in G$, and $\alpha$ is
strongly closed (open) iff, moreover, $g:X\to X$ is closed (open) on
$\domain(g)$ for each $g$.
\end{myprop}
\begin{pf} 
We prove only the closed case. Let $\domain(g)$ be closed for each $g\in
G$. We show by induction over $\length(u)$ that $\domain(u)$ is closed
for each $u\in M$: let $u$ have normal form $u=g_n\dots g_1$, so that
$\imprefix{u}=g_n\dots g_2$. Then $\domain(\imprefix{u})$ is closed by
induction. By confluence, $u\actson x$ is defined in $X$ iff $u\actson
x$ reduces to some $y\in X$. Thus, we have
\begin{displaymath}
\domain(u)=g_1^{-1}[\domain(\imprefix{u})],
\end{displaymath}
which is closed in $\domain(g_1)$ and hence in $X$. The second claim is
now trivial.
\rightqed
\end{pf}

Strong closedness is in a suitable sense `inherited' by the
globalization:

\eat{
\begin{mylem}\label{lem:closed}
If $F$ is closed (open), then the domain $\domain(w^*)$ is closed (open) in
$X$ for each $w \in W(\BC)$, and if $F$ is strongly closed (strongly
open), then $w^*:X\to X$ is a closed (open) map for each $w \in W(\BC)$.
\end{mylem}  
\begin{pf} 
We proceed by induction on $\length(w)$.  If $\length(w)=0$, then
$\domain(w^*)=X$ is closed (open). Now let $w=(g_n,\dots,g_1)$. We have
\begin{displaymath}
\domain(w^*)=g_1^{-1}[\domain((g_n,\dots,g_2)^*)]\subset \domain(g_1);
\end{displaymath}
$\domain((g_n,\dots,g_2)^*)$ is closed (open) by induction, so that
$\domain(w^*)$ is closed (open) in $\domain(g_1)$ and hence in $X$.
\rightqed
\end{pf}
}

\begin{myprop} 
\label{closed} 
\begin{rmenumerate}
\item	If $\alpha$ is strongly open, then $u: X\to Y$ is open
	for every $u \in M$.  
\item	If $\alpha$ is strongly closed and confluent, then $u: X\to Y$
	is closed for every $u \in M$. 
\end{rmenumerate}
\end{myprop}
\begin{pf} 
\emph{(i):} We have to show that $v^{-1}[u[U]]$ is open in $X$ for
each $v\in M$ and each open $U$ in $X$. We can write this set as 
\begin{displaymath}
v^{-1}[u[U]]=\bigcup_{n\in\Nat}V_{n,v},
\end{displaymath}
where $V_{n,v}$ denotes the set of all $x\in X$ such that there exists
$y\in U$ such that $(v,x)\simeq(u,y)$ is obtainable by applying the
generating relation $\sim$ of Formula~(\ref{eqn:generateeq})
(Section~\ref{preact-section}) $n$ times from left to right or from
right to left. We show by induction over $n$ that $V_{n,v}$ is open for
each $v$: the base case is trivial. Now by the definition of $\sim$,
\begin{displaymath}
V_{n+1,v}=\bigcup_{\substack{p,q\in M\\ v=pq}}
		(q^{-1}[V_{n,p}]\cap X)
	\cup
	\bigcup_{p\in M}(p[V_{n,vp}]\cap X),
\end{displaymath}
where the first part of the union corresponds to the first step in the
derivation of $(v,x)\simeq(u,y)$ being of the form $(v,x)=(pq,x)\sim
(p,q\actson x)\in V_{n,p}$ and the second to that step being of the form
$V_{n,vp}\owns (vp,z)\sim (v,p\actson z)=(v,x)$. By the inductive
assumption, the sets $V_{n,p}$ and $V_{n,vp}$ are open; hence, all
components of the union are open, since all $p\in M$ have open domains
and are open as partial maps $X\to X$.

\emph{(ii):} The argument is analogous to the one above, noticing that
thanks to confluence, all unions above can be restricted to finite ones:
the derivation of $(v,x)\simeq(u,y)$ needs at most
$\length(v)+\length(u)$ steps; in the first part of the union in the
decomposition of $V_{n+1,v}$, the decompositions $v=pq$ can be
restricted to be normal; and in the second part of the union, $p$ need
only range over generators that occur in the normal form of $u$.
\rightqed
\end{pf} 

\begin{mycor} 
\label{open} 
Let $\alpha$ be strongly open. Then the translation map $u: Y\to Y$ is
open for every $u \in M$.
\end{mycor} 

\begin{pf} 
Let $U$ be an open subset of $Y$. We have to show that $u\actson U$ 
is open. Represent this set as
\begin{displaymath}
u\actson U=\bigcup_{v\in M} uv\actson (v^{-1}[U]\cap X).
\end{displaymath}
Now observe that each component set of the union is open. Indeed, since
$v^{-1}[U]\cap X$ is open in $X$, Proposition~\ref{closed}, (i), implies
that $uv\actson (v^{-1}[A]\cap X)$ is open in $Y$.
\rightqed
\end{pf} 

As the following example shows, the `closed version' of the last statement
fails to be true even for confluent partial actions. 

\begin{rmexmp}
\label{not closed}
Let $X=\Reals$ be the real line. For $n\in \Nat$, let $p_n:
\Nat \to \Nat$ be the constant map with value $n$. These maps, together
with the identity map on $\Nat$, form a monoid $M$ which acts on
$\Nat\subset\Reals$ and thus partially acts on $X$.  Clearly, this
partial action is strongly closed; but the translation $p_1: Y \to Y$ of
the corresponding globalization fails to be closed.  Indeed, define a
subset of $Y$ as
\begin{displaymath}
A=\{p_n\actson \frac{1}{n} \mid n\ge 2\}.
\end{displaymath}
Then $A$ is closed in $Y$ because $v^{-1}[A]\cap X$ has at most one point for
every $v\in M$.  However, $p_1\actson A$ is not closed.  To see
this, observe that $p_1p_n=p_1$ and hence $p_1\actson A=\{p_1\actson
\frac{1}{n} \mid n \ge 2 \}$.  The sequence of points $p_1 \actson
\frac{1}{n}$ in $p_1\actson A$ converges to the point $p_1 \actson 0=1$,
which is outside of $p_1\actson A$.
\end{rmexmp}

\begin{rmrem}
In the case that $M$ is a group, closed partial actions are
automatically strongly closed.  Moreover, since in this case each translation
$u:Y \to Y$ is a homeomorphism, the `closed version' of
Corollary~\ref{open} is trivially true.
\end{rmrem}

We now approach the question of normality and dimension.
Let $Z$ be a topological space.  Following Wallace \cite{Wallace45}, we
say that $X$ is of {\it dimensional type} $Z$ (in short: $X\tau Z$) if,
for each closed set $A \subset X$ and each continuous map $f: A \to Z$,
there exists a continuous extension $\overline{f}: X \to Z$.

\begin{mythm} 
\label{tauthm} 
If $\alpha$ is closed and confluent, then $X \tau Z$ implies $Y \tau Z$.  
\end{mythm}
\begin{pf} 
Let $A \subset Y$ be closed, and let $\psi: A \to Z$ be a continuous
map. In order to define the required extension $\overline{\psi} : Y \to
Z$, we construct a sequence of continuous functions $\psi_n:Y_n\to Z$
(cf. Section~\ref{sec:globalizations}) such that each $\psi_n$ extends
the restriction $\restr{\psi}{A\cap Y_n}$ and each $\psi_{n+1}$ extends
$\psi_n$. We then obtain $\bar\psi$ as the union of the $\psi_n$.

$Y_0$ is just $X$. Since $A\cap X$  is closed in $X$, we can choose
$\psi_1$ as an extension of $\restr{\psi}{A\cap X}$ to $X$.

Now assume that we have constructed the sequence up to $n$. We define
auxiliary functions $\lambda_u:B_u\to Z$, where $B_u$ is closed in $X$,
for each $u\in M$ such that $\length(u)\leq n$ as follows: let $u$
have normal form $g_k\dots g_1$, and let $D$ be the (closed) domain of
$g_1$. The set $B_u$ is the union $D\cup u^{-1}[A]$ (hence closed), and
$\lambda_u$ is defined by
\begin{displaymath}
\lambda_u(x)=\left\{\begin{array}{ll}
	\psi_n(u\actson x),& \textrm{if $u\actson x\in Y_n$,
	and}\\
	\psi(u\actson x), & \textrm{if $u\actson x\in A$.}
	\end{array}\right.
\end{displaymath}
By assumption on $\psi_n$, $\lambda_u$ is well-defined. It is continuous
on $D$ and on $u^{-1}[A]$, hence continuous, since both these sets are
closed. 

Since $X\tau Z$, each $\lambda_u$ has a continuous extension
$\kappa_u:X\to Z$. We put
\begin{displaymath}
\psi_{n+1}(u\actson x)=\kappa_u(x)
\end{displaymath}
for each $u\in M$ with $\length(u)\leq n$ and each $x\in X$. Since
$\length(\globa)\leq n$ for any $\globa\in Y_{n+1}$ that admits more than one
such representation $\globa=u\actson x$, $\psi_{n+1}$ is well-defined. It is
continuous for fixed $u$, which implies overall continuity by definition
of the topology on $Y$; finally, it extends $\restr{\psi}{A\cap
Y_{n+1}}$ and $\psi_n$ by construction.
\rightqed
\end{pf}

\begin{mycor} \label{normality} 
If $\alpha$ is closed and confluent and $X$ is normal 
(and has $\dim(X) = n$), then $Y$ is normal (and has $\dim(Y)=n$).    
\end{mycor}
\begin{pf}
First note that $Y$ is a $T_1$-space by virtue of
Corollary~\ref{cor:T1-closed}.  Now use Theorem \ref{tauthm} and
well-known characterizations of normality (for $Z=[0,1]$) and dimension
(for $Z=S_n$) in terms of dimensional type.
\rightqed
\end{pf}

If $\alpha$ is not closed then we cannot in general expect  
the preservation of basic topological properties, such as for instance
$T_2$, in $Y$ (or, in fact, in any other globalizations):
\begin{rmexmp} 
Let $h: O\to O$ be an autohomeomorphism of an open subset $O$ of $X$.
Suppose that sequences $(x_n)$ and $(y_n)$ in $O$ both converge to the
same point in $X\setminus O$, and that $(h(x_n))$ and $(h(y_n))$
converge to points $c$ and $d$ in $X \setminus O$, respectively. If $X$
admits a \emph{Hausdorff} extension $ X \into Z$ such that $h$ extends
to a global map on $Z$, then $c=d$: in $Z$,
we have
\begin{displaymath}
c=\lim h(x_n)=h(\lim x_n)=h(\lim y_n)=\lim h(y_n)=d.
\end{displaymath}
It follows that $Y$ cannot be Hausdorff for any (even very good) $X$
that has such a subspace $O$ with $c$ and $d$ distinct.  As a concrete
example, take
$X=\integers\cup\{\infty,-\infty\}$, $O=\integers$, $(x_n)$ and $(y_n)$
the sequences of positive even and odd numbers, respectively, and
$h(n)=n$ if $n$ is even, $h(n)=-n$ otherwise.  (By way of contrast,
observe that, by Corollary~\ref{cor:T1-group}, the globalization
$Y$ of the group generated by $h$ \emph{is} $T_1$).
  
This example shows in particular that the abstract globalization
problem of~\cite[p. 294]{ChoiLim00} in general fails to have a
Hausdorff solution.
\end{rmexmp}

\section{Non-Expansive Partial Actions}

We will now move on from topology into the realm of metrics and
pseudometrics.

\begin{rmdefn}
A \emph{weak pseudometric} space is a pair $(X,d)$, where $d:X\times
X\to\Reals^+\cup\{\infty\}$ is a symmetric distance function that
satisfies the triangle inequality and $d(x,x)=0$ for each $x\in X$. A
\emph{pseudometric} space is a weak pseudometric space $(X,d)$ such that
$d(x,y)<\infty$ for all $x,y$.  A weak pseudometric space is called
\emph{separated} if $d(x,y)=0$ implies $x=y$. (Thus, a metric space is a
separated pseudometric space.) 

We will denote all
distance functions by $d$ (and the space $(X,d)$ just by $X$) where this
is unlikely to cause confusion. A function $f$ between weak pseudometric
spaces is called \emph{non-expansive} if $d(f(x),f(y))\leq d(x,y)$ for
all $x,y$.
\end{rmdefn}
We denote the categories of weak pseudometric, pseudometric, and metric spaces
with non-expansive maps as morphisms by $\PMet$,
$\SMet$, and $\Met$, respectively.

A partial action of a monoid $M$ on a weak pseudometric space $X$ is called
\emph{non-expansive} if the partial map $u:X\to X$ is non-expansive on
its domain (as a subspace of $X$) for each $u\in M$. Note here that both
$\SMet$ and $\Met$ are closed under subspaces in
$\PMet$. 
\eat{Correspondingly, a \emph{non-expansive preaction} (cf.\
Example~\ref{exmp:preact}) of a category $\BC$ on $X$ is a functor
$F:\BC\to\BM(X)$, where $\BM(X)$ denotes the category of weak pseudometric
subspaces of $X$ and non-expansive maps.
}

Since $\PMet$ is a topological category \cite{AdamekHerrlich90},
globalizations can be constructed in the same way as for topological
partial actions by means of final lifts: in general, given weak
pseudometric spaces $Y_i$, $i\in I$, and a family of maps $f_i:Y_i\to X$
into some set $X$, the final lift of $\CS=(Y_i,f_i)_I$ is the largest
weak pseudometric on $X$ (w.r.t.\ the pointwise order on real-valued
functions) that makes all the $f_i$ non-expansive maps. Explicitly,
given points $x$ and $y$ in $X$, an \emph{$\CS$-path} $\pi$ from $x$ to
$y$ of \emph{length} $n$ is a sequence
$((i_1,x_1,y_1),\dots,(i_n,x_n,y_n))$, $n\geq 1$, such that $x_j,y_j\in
Y_{i_j}$, $j=1,\dots,n$, $f_{i_1}(x_1)=x$,
$f_{i_j}(y_j)=f_{i_{j+1}}(x_{j+1})$ for $j=1,\dots,n-1$, and
$f_{i_n}(y_n)=y$. The associated \emph{path length} is
\begin{displaymath}
\sum_{j=1}^{n} d_j(x_j,y_j).
\end{displaymath}
In case $x\neq y$, the distance of $x$ and $y$ is easily seen to be
given as the infimum of the path length, taken over all $\CS$-paths from
$x$ to $y$ (in particular, the distance is $\infty$ if there is no such
path); otherwise the distance is, of course, $0$. If the $f_i$ are jointly
surjective (which they are in the case we are interested in), then there
is always a trivial $\CS$-path from $x$ to $x$, so that the
case $x=y$ does not need special treatment.  Due to the triangle
inequality, it suffices to consider paths $((i_j,x_j,y_j))$ where
$(i_j,y_j)$ is always different from $(i_{j+1},x_{j+1})$.

Now given a partial action $\alpha$ on a weak pseudometric space $X$, we
construct the underlying set of the free globalization $Y$ as in
Section~\ref{sec:globalizations} (as for topological spaces, we shall
keep the notation $\alpha$, $X$, $Y$ etc. throughout). It is easy to see
that free globalizations of partial actions on weak pseudometric spaces
(i.e.\ reflections into the full subcategory spanned by the total
actions in the category of partial actions) are, as in the topological
case, given as final lifts of the family $\CS$ of maps
\begin{displaymath}
u:X\to Y,
\end{displaymath}
where $u$ ranges over $M$. For the sake of clarity, we denote the
distance function on $Y$ thus defined by $D$.

\emph{For the remainder of this section, we shall assume that $\alpha$
is confluent.}  

Under this condition, one may further restrict the
paths to be taken into consideration: in general, we may write an
$\CS$-path $\pi$ from $\globa$ to $\globb$ ($\globa,\globb\in Y$) in the
form
\begin{displaymath}
u_1\actson x_1,u_1\actson y_1\redeq u_2\actson x_2,\dots,
u_{n-1}\actson y_{n-1}\redeq u_n\actson x_n, u_n\actson y_n
\end{displaymath}
(in short: $(u_j,x_j,y_j)$), where $u_1\actson x_1=\globa$ and
$u_n\actson y_n=\globb$ .  Denote by $D(\pi)$ the corresponding path
length $\sum_{j=1}^{n} d(x_j,y_j)$.  By definition, $D(\globa, \globb) =
\inf D(\pi) $ where $\pi$ runs over all possible paths. Recall that
 $D(\globa, \globb) =\infty$ iff there is no path from
$\globa$ to $\globb$.  We say that $\pi$ is {\it geodesic} if
$D(\globa, \globb) = D(\pi)$.

There are two additional \emph{assumptions} we may introduce: 
\begin{rmenumerate}
\setlength{\itemsep}{1em}
\item	\emph{For each $j=1,\dots,n$, at least
	one of $u_j\actson x_j$ and $u_j\actson y_j$ is in normal form.} 
	
	~

        \noindent Indeed, if $u_j$ has normal form $g_k\dots g_1$ and
	both $x_j$ and $y_j$ are in the domain of $g_1$, then we obtain
	a shorter path replacing $(u_j,x_j,y_j)$ by $(g_k\dots
	g_2,g_1\actson x_j,g_1\actson y_j)$ (since $g_1$ is
	non-expanding).
\item	\emph{For each $j=1,\dots,n-1$, at most one of $u_j\actson y_j$ and 
	$u_{j+1}\actson x_{j+1}$ is normal.}

	~

        \noindent By the above, we may assume $(u_j,y_j)\neq(u_{j+1},x_{j+1})$. 
	But both these pairs represent the same 
        point of $Y$, which has only one normal form.
\end{rmenumerate} 
We will henceforth consider only $\CS$-paths that are \emph{reduced}
according to these assumptions. We denote the transitive closure of the
one-step reduction $\redstep$ by $\redgt$ (reversely: $\redlt$); i.e.\
$\redgt$ is like $\redge$ except that we require that at least one
reduction step takes place. If $u_j\actson y_j$ is reducible and
$u_{j+1} \actson x_{j+1}$ is normal then necessarily $u_j\actson y_j
\redgt u_{j+1}\actson x_{j+1}$, which we will indicate in the notation
for paths; similarly if $u_j\actson y_j$ is normal and $u_{j+1}\actson
x_{j+1}$ is reducible.

The `normality patterns' that occur in reduced
paths are restricted in a rather amusing way:

\begin{mylem}\label{lem:reduced-path}
Every reduced path from $\globa \in Y$ to $\globb \in Y$  
has one of the following forms:
\begin{enumerate}
\renewcommand{\labelenumi}{{\rm (A\arabic{enumi})}}
\setlength{\labelsep}{0cm}
\item  \hskip 0.5cm 
$n, r \redgt \cdots \redgt n, r$ 
\item  \hskip 0.5cm 
$r,n \redlt \cdots \redlt r,n$  
\item  \hskip 0.5cm
$n, n$
\item    \hskip 0.5cm 
$n, r \redgt \cdots \redgt n, r \redgt n, n$ 
\item   \hskip 0.5cm 
$n,n \redlt r,n \redlt \cdots \redlt r,n$ 
\item  \hskip 0.5cm 
$n, r \redgt \cdots \redgt n,r \redgt n, n \redlt r,n \redlt \cdots \redlt
r,n$ 
\item  \hskip 0.5cm
$n, r \redgt \cdots\redgt n,r \redeq r, n \redlt \cdots \redlt
r,n$
\end{enumerate}
where `n' and `r' mean that the corresponding term 
of the  path is normal or reducible, respectively. (Patterns such as
$n,r\redgt\cdots\redgt n,r$ are to be understood as `one or more
occurences of $n,r$'.)
\end{mylem}
\begin{pf}
If the path does not contain either of the patterns $n,n$ and $r\redeq
r$, then it must be of one of the forms $(A1)$ and $(A2)$. The occurence
of $n,n$ in some place determines the entire pattern due to restrictions
(i) and (ii) above, so that the path has one of the forms
$(A3)$--$(A6)$. Similarly, a path that contains the pattern $r\redeq r$
must be of the form $(A7)$.
\eat{Take into account that if some $u_j\actson x_j$ is reducible 
(= ``r'') then the remaining right part of this path 
necessarily looks as an  
alternating sequence of normal and reducible terms. 
That is, $u_k\actson x_k$ is reducible and
$u_k\actson y_k$ is normal for $k\geq j$. 
Indeed,   
by observation (i), $u_j\actson y_j$ must be normal. 
Now, by observation (ii), $u_{j+1}\actson x_{j+1}$ is reducible and 
etc.} 
\rightqed
\end{pf}

\eat{
\begin{mylem}\label{lem:XY-reduced-path}
Let $\pi=((u_j,x_j,y_j))$ be a reduced path of length $n$ from $x\in X$
to $y\in Y- X$. Then $\pi$ is of type $(A3)$ or $(A5)$. Therefore, 
$u_j\actson y_j$ is normal for $j=1,\dots,n$, and
$u_j\actson x_j$ is reducible for $j=2,\dots,n$.
\end{mylem}
\begin{pf}
Observe that normality of $u_1\actson x_1$ implies $u_1=e$ and hence
normality of $u_1\actson y_1$, so that $u_1\actson y_1$ must be normal
in any case. This implies that $u_j\actson x_j$ is reducible and
$u_j\actson y_j$ is normal for $j\geq 2$. 
If $n=1$, then $u_1\neq e$; i.e.\ $u_1\actson x_1$ is reducible (this
covers the remaining case $j=n=1$).
\end{pf}
} 

A first consequence of this lemma is that every space is a subspace
of its globalization:

\begin{mylem}\label{lem:XXpath}
Let $x,y\in X$. Then $((e,x,y))$ is the only reduced path from $x$ to
$y$. 
\end{mylem}
\begin{pf}
Since $e\actson z$ is in normal form for all $z\in X$,
any reduced path from $x$ to $y$ must have form $(A3)$ of
Lemma~\ref{lem:reduced-path} (all other forms either begin with
the pattern $n,r$ or end with $r,n$).
\rightqed
\end{pf}

\begin{mythm}\label{thm:pseudo-isometric}
The embedding $X\into Y$ of a weak pseudometric space into its free
globalization is isometric.
\end{mythm}
\begin{pf} 
Immediate from Lemma~\ref{lem:XXpath}.
\rightqed
\end{pf}

Of course, we are mainly interested in \emph{metric} globalizations.
Now any weak pseudometric space has a separated reflection obtained by
identifying points with distance zero. If $X$ is a separated space, then
the separated reflection $\bar Y$ of $Y$ is the free \emph{separated}
globalization of $X$, and $X$ is isometrically embedded in $\bar Y$,
since its points have positive distances in $Y$ and are hence kept
distinct in $\bar Y$. We will see below (Theorem~\ref{thm:positive})
that working with the separated reflection is unnecessary for closed
partial actions. Finiteness of distances is, on the one hand, more
problematic since there is no universal way to transform a weak
pseudometric space into a pseudometric space. On the other hand,
finiteness of distances is preserved in most cases:
\begin{rmdefn}\label{def:degenerate}
$\alpha$ is called \emph{nowhere degenerate}
if $\domain(g)\neq\emptyset$ for each $g\in G$.
\end{rmdefn}
\begin{myprop}\label{prop:pseudometric}
If $X$ is a non-empty pseudometric space, then $Y$ is pseudometric iff $\alpha$
is nowhere degenerate.
\end{myprop}
\begin{pf}
If $\alpha$ is nowhere degenerate, then there exists, for each $y\in Y$,
a path from $y$ to some $x\in X$; hence, there is a path between any two
points of $Y$, so that the infimum defining the distance function on $Y$
is never taken over the empty set and hence never infinite. If,
conversely, $\domain(g)=\emptyset$ for some $g\in G$, then there is no
reduced path (and hence no path at all) from $x$ to $g\actson x$ for
$x\in X$, so that $D(x,g\actson x)=\infty$. Indeed, assume that $\pi$ is
such a path. Since both $e\actson y$ and, by assumption on $g$,
$g\actson y$ are normal for all $y\in X$, the normality pattern of $\pi$
as in Lemma~\ref{lem:reduced-path} can neither begin with $n,r$ nor end
with $r,n$. Thus, $\pi$ must be of the form $(A3)$, which is impossible
since $\domain(g)=\emptyset$ implies $g\neq e$.
\rightqed
\end{pf}

\begin{rmrem}
Another approach to the problem of infinite distances is to consider
only spaces of diameter at most $1$ and put $D(x,y)=1$ for $x,y\in Y$ in
case there is no path from $x$ to $y$. 
\end{rmrem}

\begin{rmobserve}
Let $\globa,\globb \in Y$ have normal forms $\globa=u\actson x$ and
$\globb=v\actson y$, and let $\pi$ be a reduced path from $\globa$ to
$\globb$. If $\pi$ is of the form (A2) or (A5) of
Lemma~\ref{lem:reduced-path}, then necessarily $u \prefixle v$, and if
$\pi$ is of the form (A1) or (A4), then $v\prefixle u$. Clearly, if
$\pi$ is of the form (A3) then $u=v$.  Thus, if $u$ and $v$ are
incomparable under $\prefixle$ then $\pi$ must be of the form (A6) or
(A7).
\end{rmobserve}

\begin{mylem} \label{lem:distance} 
Let $\globa,\globb\in Y$ have normal forms $\globa=u\actson x$ and
$\globb= v\actson y$, where $u$ has normal form $g_k\dots g_1$.
\begin{rmenumerate}
\item	 If $D(\globa, \globb) < d(x, \domain(g_1))$,
	then  $u \prefixle v$.  
\item	If $u=v$ then  
	\begin{displaymath}
	\min\{d(x,y),d(x,\domain(g_1))+d(y,\domain(g_1))\}\le
	D(\globa,\globb)\le
	d(x,y).
	\end{displaymath} 
\end{rmenumerate}
\end{mylem} 
\begin{pf} Let $\pi$ be a reduced path from $\globa$ to $\globb$.

\emph{(i):} $\pi$ cannot have a normality pattern of the form
$n,r\redgt\dots$, since in that case, the first step of the path would
already contribute at least $d(x,\domain(g_1))$ to
$D(\globa,\globb)$. Hence, $\pi$ must be of one of the forms (A2), (A3),
or (A5) of Lemma~\ref{lem:reduced-path}. By the observation above, this
implies $u\prefixle v$. 

\emph{(ii):} $\pi$ must have one of the forms (A3), (A6), or (A7) of
Lemma~\ref{lem:reduced-path}. In the case (A3), $D(\pi)=d(x,y)$.  In the
cases (A6) and (A7), the normality pattern of $\pi$ is of the form $n,r
\redgt \dots \redlt r,n$. Therefore $D(\pi)\ge
d(x,\domain(g_1))+d(y,\domain(g_1))$. This proves the first inequality;
the second follows from the fact that $u:X\to Y$ is non-expansive.
\rightqed
\end{pf}

We say that a function $\phi: E \to L$ between pseudometric spaces is
\emph{locally isometric} if for every $x\in E$ there exists $\varepsilon >0$
such that $\phi$ isometrically maps the $\varepsilon$-ball
$B(x,\varepsilon)$ in $E$ onto the $\varepsilon$-ball
$B(\phi(x),\varepsilon)$ in $L$.  Clearly, $E$ is separated
iff $\phi(E)$ is separated.  Every locally isometric injective map is a
topological embedding.

\begin{myprop}\label{prop:wc-layers} 
If $\alpha$ is closed, then
\begin{rmenumerate}
\item	$D(u \actson x, v \actson y)=0$ implies	$u=v$ for normal 
	forms $u\actson x$, $v\actson y$.
\item	The set $\bigcup_{u\prefixle v} v\actson R_v$ is open for each $u$.
\item	Each $Y_k$ (in particular, $Y_0=X$) is closed in $Y$. 
\item	The subspace $Y_{k+1}\setminus Y_k$ is a topological sum
	$\bigcup_{\length(u)=k+1} u\actson R_u$ of disjoint subsets 
	$u\actson R_u$. 
\item	For every $u\in M$ the bijective function $u: R_u \to u\actson R_u$ is
	locally isometric (and, hence, a homeomorphism). 
\end{rmenumerate}
\end{myprop} 
\begin{pf} 
\emph{(i):} Let $u$ have normal form $g_n\dots g_1$. Then $D(u\actson
x,v\actson y)=0<d(x,\domain g_1)$ by closedness, so that $u\prefixle v$ by
Lemma~\ref{lem:distance}~(i). Analogously, $v\prefixle u$.

\emph{(ii):} Let $u\in M$, and let $a$ have normal form $p\actson x$
(i.e.\ $a\in p\actson R_p$) for some $u\prefixle p$ with normal form
$p=g_n\dots g_1$. Put $\eps=d(x,\domain(g_1))$. By closedness,
$\eps>0$. By Lemma~\ref{lem:distance}~(i), the $\eps$-neighbourhood of
$a$ is contained in $\bigcup_{p\prefixle v}v\actson R_v$ and hence in
$\bigcup_{u\prefixle v}v\actson R_v$, which proves the latter set to be
open.

\emph{(iii):} The complement of $Y_k$ is a union of sets
$\bigcup_{u\prefixle v} vR_v$. 

\emph{(iv):} Disjointness is clear, and by (ii), each set $u\actson R_u$
with $\length(u)=k+1$ is open in $Y_{k+1}\setminus Y_k$, since $u\actson
R_u=\left(\bigcup_{u\prefixle v}v\actson R_v\right)\cap(Y_{k+1}\setminus Y_k)$.

\emph{(v):} Let $u=g_k \dots g_1$ be normal, and let $x \in
R_u=X\setminus\domain(g_1)$. Since $\alpha$ is closed,
$\varepsilon:=d(x,\domain(g_1))>0$. By Lemma \ref{lem:distance}~(ii),
the bijective function $u: R_u \to u\actson R_u$ isometrically maps the
$\varepsilon$-ball $B(x,\varepsilon)$ onto the $\varepsilon$-ball
$B(u\actson x,\varepsilon)$ in $u\actson R_u$.
\rightqed
\end{pf}

As an immediate consequence, we obtain the announced separatedness result:

\begin{mythm}\label{thm:positive}
If $\alpha$ is closed and $X$ is separated, then $Y$ is separated.
\end{mythm} 
\begin{pf} 
Let $u\actson x$ and $v\actson y$ be normal forms in $Y$ with
$D(u\actson x,v\actson y)=0$. Then $u=v$ by
Proposition~\ref{prop:wc-layers}~(i);  therefore $x,y\in R_u$. By
Proposition~\ref{prop:wc-layers}~(v), $D(u\actson x,u\actson y)=0$
implies $d(x,y)=0$ and hence $x=y$.
\rightqed
\end{pf} 

\eat{
Assume $m\neq n$. Let $m$ and $n$ have normal forms $g_k\dots g_1$ and
$g_r\dots g_1$, respectively. 
Let $\pi=((u_j,x_j,y_j))$ be a reduced 
path of length $l$ from $m\actson x$ to $n\actson y$.   
First suppose that $\pi$ is of type $(A2), (A5)$ or $(A6)$. Then 
$u_l\actson x_l$ is reducible and $u_l\actson y_l$ is normal, i.e.\
$e\neq u_l=n$. Since $F$ is closed, 
$y$ has positive distance $\epsilon$ from the domain $\domain(g_1)$ of $g_1$. 
Therefore, $\pi$ has path length at least $\epsilon$.   
Similarly, if $\pi$ is of type $(A1)$ or $(A4)$, then $\pi$ has path
length at least $\delta$, where $\delta>0$ is the distance of $x$ from the
domain of $g_1$. Type $(A3)$ is excluded since $m\neq n$. Thus, 
\begin{displaymath}
d(m\actson x,n\actson y)\geq\min \{\epsilon,\delta \} >0.
\end{displaymath}
}

\eat{
Let $m$ have normal form $g_k\dots g_1$. Since $X$ is separated,
we can, by Theorem~\ref{thm:pseudo-isometric}, assume $k>0$ . Let $\pi$
be a reduced path from $m\actson x$ to $m\actson y$. If $\pi$ has the
form $(A3)$ of Lemma~\ref{lem:reduced-path}, then its path length is
$d(x,y)>0$. All other forms either begin with $n,r$ or end with $r,n$,
so that in these cases, the path length of $\pi$ is at least the
positive distance of the set $\{x,y\}$ from the domain of $g_1$.
\end{pf}
}

\begin{rmrem}
The converse of the above theorem holds if $X$ is complete: assume that
$Y$ is separated, let $g\in G$, and let $(x_n)$ be a convergent sequence
in $\domain(g)$; we have to show that $x=\lim x_n$ is in $\domain(g)$.
Now $(g\actson x_n)$ is a Cauchy sequence in $X$, hence by assumption
convergent; let $z=\lim g\actson x_n$. For every $n$, we have a path
\begin{displaymath}
e\actson z, e\actson (g\actson x_n) \redlt g\actson x_n,
g\actson x
\end{displaymath}
from $z$ to $g\actson x$. The associated path length is $d(z, g(x_n)) +
d(x_n, x)$, which converges to $0$ as $n\to\infty$. Hence, $D(z,
g\actson x)=0$, so that $z=g\actson x$ by separatedness; this implies
that $g\actson x$ is defined in $X$ as required.
\end{rmrem}

\begin{rmexmp}
Even for closed partial actions of groupoids on metric spaces, the
metric globalization does not in general induce the topology of the
topological globalization of Section~\ref{sec:globalizations}. Take, for
instance, $X=[0,1]$. The full subcategory of $\BM(X)$ spanned by all
singleton subspaces induces a partial action $\alpha$ as described in
Example~\ref{exmp:confluentpa}~\romanref{exmp:preact} (cf.\ also
Example~\ref{exmp:shimrat}). The universal \emph{topological}
globalization $Y$ of $\alpha$ is not even first countable: as in
Example~\ref{exmp:shimrat}, denote the map $\{x\}\to\{y\}$ by $(yx)$ for
$x\neq y$ in $X$. Then we have a subspace $Z$ of $Y$ formed by all
points of the form $x$ or $(y0)\actson x$. The space $Z$ is
homeomorphic to the quotient space obtained by taking one base copy of
$[0,1]$ and uncountably many copies of $[0,1]$ indexed over the base
copy, and then identifing for each $a\in[0,1]$ the $0$ in the $a$-th
copy with the point $a$ in the base copy. In particular, already $Z$
fails to be first countable.
\end{rmexmp}

\begin{mythm} \label{thm:metricdim} 
If $X$ is a metric space and $\alpha$ is closed and nowhere degenerate, then
$Y$ is a metric space. Moreover, $\dim(Y)=\dim(X)$.  
\end{mythm} 
\begin{pf}
By Theorem~\ref{thm:positive} and
Proposition~\ref{prop:pseudometric}, $Y$ is a metric space. 

It remains to be shown that $\dim(X)=\dim(Y)$. Now $Y=\cup_{n \in \Nat}
Y_n$ where, by Proposition~\ref{prop:wc-layers}, each $Y_n$ is a closed
subset of $Y$. Therefore, by the standard \emph{countable sum theorem},
it suffices to show that $\dim(Y_n) \leq \dim(X)$ for every $n$. We
proceed by induction.  The case $n=0$ is trivial, since $Y_0=X$. We have
to show that $\dim(Y_{n+1}) \leq \dim(X)$ provided that $\dim(Y_n) \leq
\dim(X)$.  The idea is to use the following result of Dowker~\cite{Dowker55}.

\begin{mylem}[Dowker]  
Let $Z$ be a normal space, and let $Q$ be a closed 
subspace of $Z$ such that 
$\dim(Q)\leq k$. Then $\dim(Z) \leq k$ 
if and only if every closed subspace 
$A\subset Z$ disjoint from $Q$ satisfies $\dim(A) \leq k$. 
\end{mylem}

We apply this lemma to the closed subspace $Y_n$ of $Y_{n+1}$. By the
induction hypothesis, $\dim(Y_n) \leq \dim(X)$.  We have to show that
$\dim(A) \leq \dim(X)$ for every closed subset $A$ of $Y_{n+1}$ which is
disjoint from $Y_n$, i.e.\ $A \subseteq Y_{n+1} \setminus Y_n$. By
Proposition \ref{prop:wc-layers} (iv), $A$ is a topological sum
$\bigcup_{\length(u)=n+1} A_u$ of disjoint subspaces $A_u:=A \cap
u\actson R_u$. Each $A_u$ is a subspace of $u\actson R_u$. Therefore, by
Proposition~\ref{prop:wc-layers}~(v), $A_u$ is homeomorphic to a
subspace of $X$. Since the dimension is hereditary (for arbitrary, not
necessarily closed subspaces) in perfectly normal (e.g., metrizable)
spaces, we have $\dim(A_u) \leq \dim(X)$. Thus, $\dim(A) \leq
\dim(X)$. By Dowker's result this yields $\dim(Y_{n+1}) \leq \dim(X)$.
\rightqed
\end{pf}

\begin{rmrem}
One application of Theorems~\ref{thm:pseudo-isometric}
and~\ref{thm:metricdim} is to obtain all sorts of metric gluing
constructions. A simple example of this is Theorem~2.1
of~\cite{Bogatyi02}, which states that \emph{given metric spaces $X_1$
and $X_2$ with intersection $Z=X_1\cap X_2$ such that $Z$ is closed both
in $X_1$ and in $X_2$ and the metrics of $X_1$ and $X_2$ agree on $Z$,
there exists a metric on $X_1\cup X_2$ which agrees with the given
metrics on $X_1$ and $X_2$, respectively}. Using our results, this can
be seen as follows: let $G$ be the free group with a single generator
$u$ (i.e.\ $G\cong\integers$), let $X_1+X_2=X_1\times\{1\}\cup
X_2\times\{2\}$ be the disjoint union of $X_1$ and $X_2$, and let a
partial action of $G$ on $X_1+X_2$ be defined by $u\actson(x,1)=(x,2)$
(and $u^{-1}\actson(x,2)=(x,1)$) for $x\in Z$. This partial action is
closed and, by
Example~\ref{exmp:confluentpa}~\romanref{exmp:freegrouppa},
confluent. In the globalization $Y$, we find the set $X_1\cup X_2$
represented as $W=(u\actson X_1)\cup X_2$, and the metric on $W$ agrees
with the respective metrics on $X_1$ and $X_2$, since the maps
$f_1:X_1\to W$ and $f_2:X_2\to W$ defined by $f_1(x)=u\actson (x,1)$ and
$f_2(y)=(y,2)$ are isometries.
\end{rmrem}


In standard terminology, some of the above results can be summed
up as follows:
\begin{mythm} 
\label{mglobal}
Let $\Gamma$ be a set of partial non-expansive maps (isometries) with
non-empty closed domain of a metric space $X$. Then there exists a
closed isometric embedding $X \hookrightarrow Y$ into a metric space $Y$
such that all members of $\Gamma$
can be extended to global non-expansive maps (isometries) of $Y$ and such
that, moreover, $\dim(Y)=\dim(X)$ and $|Y| \leq |X| \actson | \Gamma | 
\actson \aleph_0$.
\end{mythm}
\begin{pf} 
$\Gamma$ generates a subcategory (a subgroupoid, if all members of
$\Gamma$ are partial isometries) $\BC$ of the category $\BM (X)$ of
metric subspaces of $X$; the set of morphisms of $\BC$ has cardinality
at most $ | \Gamma | \actson \aleph_0$.  The inclusion $\BC
\hookrightarrow \BM (X)$ induces a closed non-expansive nowhere
degenerate partial action $\alpha$ on $X$ as described in
Example~\ref{exmp:confluentpa}~\romanref{exmp:preact}. By Theorem
\ref{thm:metricdim} and Proposition \ref{prop:wc-layers}~(iii), the
globalization of $X$ w.r.t.\ $\alpha$ has the desired properties.
\rightqed
\end{pf}

By iterating the construction above, we can improve, in part, the known
result \cite{Uspenskij98}\footnote{Uspenskij shows that it can be
assumed that the weight is preserved and that the isometry group of $X$
(endowed with the pointwise topology) is topologically embedded into the
isometry group of~$Z$ (but this construction does not preserve dimension)}
that every metric space $X$ can be
embedded into a metrically ultrahomogeneous 
space $Z$:
\begin{mythm} 
\label{thm:sigma-homog} 
For every metric space $X$ there exists an isometric closed embedding
$X\hookrightarrow Z$ into a metrically ultrahomogeneous space $Z$
such that $\dim(Z) = \dim(X)$ and $|Z|=|X|$. 
\end{mythm}
\begin{pf}
Start with the set $\Gamma$ containing all partial isometries between
finite subspaces of $X$ \emph{and} all global isometries of $X$ (here,
$\Gamma$ is already a subcategory of $\BM(X)$). Let $Z_1$ be the
corresponding globalization according to the above theorem
and iterate this process; the direct limit $Z_{\infty}$ of the resulting
ascending chain of metric spaces $X\hookrightarrow Z_1 \hookrightarrow
Z_2 \hookrightarrow \dots $ is an ultrahomogeneous space. Moreover, each
inclusion is closed and $\dim(Z_n)=\dim(X)$ for all $n$. Hence, the
inclusion $X\into Z_\infty$ is closed, and by the countable sum theorem,
$\dim(Z_\infty)=\dim(X)$. A more careful choice of global isometries
will guarantee that $|Z|=|X|$.
\rightqed
\end{pf}

\begin{rmrem} 
Topological versions of Theorems~\ref{mglobal}
and~\ref{thm:sigma-homog}, with `metric' replaced by `normal' and
`metrically ultrahomogeneous' by `topologically ultrahomogeneous', can
be derived using Corollary~\ref{normality} (see also~\cite{Megrelishvili85,Megrelishvili86,Megrelishvili00}).
\end{rmrem} 

The global metric $D$ on $Y$ is in 
some respects easier to handle in case $M$ is a group.
Since the
elements of $M$ act as isometries and hence $D(u\actson x, v
\actson y)=D(x,u^{-1}v\actson y)$ for all $u,v\in M$ and all
$x,y\in X$, it suffices to consider distances of the form $D(x,u\actson
y)$. Thus, the calculation of distances can be simplified:

\begin{myprop}\label{groupoid} 
Let $M$ be a group. Let $u, v\in M$, let $g_k\dots g_1$ be
the normal form of $u^{-1}v$, and let $x,y\in X$. Then
\begin{displaymath}
D(u\actson x, v \actson y)=
\inf \left( d(y,x_1)+\sum_{i=1}^k d(g_i(x_i), x_{i+1}) \right), 
\end{displaymath}
where $x_i$ ranges over $\domain(g_i)$ for $i =1,\dots,k$ and
$x_{k+1}=x$.  
\end{myprop} 
\begin{pf}  
As explained above, we need only calculate the distance from
$\globa:=u^{-1}v\actson y$ to the point $x\in X$. 

Since $e\actson z$ is normal for all $z\in X$, a reduced path $\pi$ from
$\globa$ to $x$ cannot end with the normality pattern $r,n$, so that
(excluding the trivial case (A3)) $\pi$ must have one of the forms
$(A1)$ or $(A4)$ of Lemma~\ref{lem:reduced-path}. Thus, $\pi$ is
determined by a subdivision $s_r\dots s_1$ of $(g_k,\dots,g_1)$ into
non-empty words $s_i$ and a selection of elements $x_i\in \domain(s_i^*)$,
$i=1,\dots,r$; putting $x_{r+1}=x$, we can write the corresponding path
length as
\begin{displaymath}
d(y, x_1)+\sum_{i=1}^r d(s_i^*(x_i),x_{i+1}).
\end{displaymath} 
Now observe that one subdivision of $(g_k,\dots,g_1)$
is that into $k$ one-element subwords $s_i=(g_i)$. Selecting elements
$x_i\in\domain(s_i^*)=\domain(g_i)$, $i=1,\dots,k$, defines a (not
necessarily reduced) path; call such paths \emph{elementary paths}. It
is easy to see that any reduced path $\pi$ gives rise to an elementary
path $\bar\pi$ such that $D(\pi)=D(\bar\pi)$, and the lengths of
elementary paths are exactly the sums given in the formula of the
statement.
\rightqed
\end{pf}

A further rather immediate consequence of Lemma~\ref{lem:reduced-path}
is the existence of geodesic paths under suitable compactness
assumptions: 

\begin{rmdefn}
Let $u\in M$ have normal form $g_k\dots g_1$, $k\ge 0$.  $u$ is
called a \emph{$C$-element} if $\domain(g_i)$ is 
compact for $i=1,\dots,k$.  A partial action is \emph{compact} if
$\domain(f)$ is compact for every morphism $f$.
\end{rmdefn}
Clearly, $\alpha$ is compact iff every $u \in M$ is a $C$-element.
\begin{mythm}\label{geodesic} 
Let $X$ be a 
weak pseudometric space. If $u$ and $v$ are
$C$-elements and $\globa=u\actson x$, $\globb=v\actson y$ are normal,
then there exists a geodesic from $\globa$ to $\globb$. In particular, if
$\alpha$ is compact then there exists a geodesic for every pair of elements
in $Y$.
\end{mythm}
\begin{pf} 
It suffices to show that, for each of the forms listed
in Lemma~\ref{lem:reduced-path}, there exists a path which realizes the
infimum among all reduced paths of that form. We treat only the case
(A7); the other cases are analogous (and, mostly, easier). 

A reduced
path $((u_j,x_j,y_j))$ from $\globa$ to $\globb$ of the form (A7) is
determined by a choice of a sequence $(u_1,\dots,u_k)$ such that
\begin{displaymath}
u=u_1\prefixgt \dots \prefixgt u_r\quad\textrm{and}\quad 
u_{r+1}\prefixlt \dots \prefixlt u_k=v
\end{displaymath}
for some $1\le r \le k-1$, and a choice of elements $y_i\in \domain(g^i_1)$,
$i=1,\dots,r$ and $x_i\in \domain(g^i_1)$, $i=r+1,\dots,k$, where $u_i$ has
normal form $g^i_{s_i}\dots g^i_1$. Obviously, there are only finitely
many choices of $(u_1,\dots,u_k)$, so that it suffices to show that,
given such a choice, the infimum among the corresponding paths is
realized by some choice of elements as described. This follows by a
standard compactness argument: the $\domain(g^i_1)$ are compact, and the
path length depends continuously on the choice of the $x_i$ and $y_i$.
\rightqed
\end{pf}

\begin{mycor}
Let $\alpha$ be compact. If $X$ is a \emph{path space}, i.e.\ if the
distance between any two points is the infimum of the lengths of all
curves joining the points~\cite{Gromov99}, then so is $Y$.
\end{mycor}

\section{Conclusion and Outlook}

We have demonstrated how a simple set-theoretic construction of
globalizations for partial actions of monoids can be applied to
topological and metric spaces, and we have shown that the resulting
extensions are surprisingly well-behaved, provided that the partial
action is confluent. In particular, we have shown
that, in both cases, the original space is embedded in its extension,
and that, under natural assumptions, important properties such as
dimension, normality, or path metricity are preserved. Classical
homogenization results arise as special cases of our construction.
The main tool has been the application of rewriting theory in order to gain
better control of the globalization.

Open questions include preservation of further topological and metric
properties by the globalization, as well as the extension of the method
to other categories. This includes categories used in general topology
such as uniform spaces or, more generally, nearness
spaces~\cite{Herrlich74}, as well as, in the realm of distance
functions, the category of approach spaces~\cite{Lowen97}, but also
structures of a more analytical nature such as measurable maps (of
$mm$-spaces \cite{Gromov99,Pestov02}), smooth maps, or conformal maps.

\begin{ack}
We wish to thank Horst Herrlich and the anonymous referee for valuable
suggestions.
\end{ack}

\bibliographystyle{amsplain}
\bibliography{partact}

\end{document}

%% file: header.tex
\newcounter{tempref}
\newcommand{\romanref}[1]{\setcounter{tempref}{0\ref{#1}}(\roman{tempref})}

\newenvironment{rmenumerate}%
	{\begin{enumerate}
	\renewcommand{\labelenumi}{(\roman{enumi})}} %
	{\end{enumerate}}
\newenvironment{alenumerate}%
	{\begin{enumerate}\renewcommand{\labelenumi}{(\alph{enumi})}}%
	{\end{enumerate}}
	{\begin{enumerate}\renewcommand{\labelenumi}{\arabic{enumi})}}%
	{\end{enumerate}}
\newenvironment{mypf}[1]%
  {\begin{pf*}{PROOF (#1).}}{\end{pf*}}
	{\setlength{\itemsep}{0pt}%
	\setlength{\topsep}{0pt}%
	\begin{trivlist}\item[]\begin{displaymath}%
	\begin{array}{rcll}}%
	{\end{array}\end{displaymath}\end{trivlist}}
	{\begin{list}{}{\setlength{\leftmargin}{0pt}%
			\setlength{\topsep}{1cm}}%
	\item{\bf Acknowledgements:}}%
	{\end{list}}

\newcommand{\rightqed} {\mystrut{0em}\hfill{$\Box$}}

\newcommand{\Cat}{\mathbf}
\newcommand{\Cls}{\mathcal}
\newcommand{\Op}{\mathrm}

\newcommand{\Mor}{{\Op{Mor}\,}}

\newcommand{\BB}{{\Cat B}}
\newcommand{\BC}{{\Cat C}}

\newcommand{\BM}{{\Cat M}}
\newcommand{\BS}{{\Cat S}}
\newcommand{\BT}{{\Cat T}}

\newcommand{\PMet}{{\Cat{wPMet}}}
\newcommand{\SMet}{{\Cat{PMet}}}
\newcommand{\Met}{{\Cat{Met}}}

\newcommand{\CS}{{\Cls S}}

\newcommand{\eps}{\varepsilon}

\newcommand{\integers}{{\mathbb Z}}

\newcommand{\openbox}{\leavevmode
  \hbox to.77778em{%
  \hfil\vrule
  \vbox to.675em{\hrule width.6em\vfil\hrule}%
  \vrule\hfil}}
\newcommand{\mystrut}[1]{\rule[#1]{0cm}{0.1cm}}

\newcommand{\into}{\hookrightarrow}
\newcommand{\id}{{id}}

\newcommand{\domain}{\operatorname{dom}}

\newcommand{\actson}{\cdot}
\newcommand{\length}{\operatorname{lg}}
\newcommand{\Nat}{{\mathbb N}}
\newcommand{\Rat}{{\mathbb Q}}
\newcommand{\Reals}{{\mathbb R}}

\newcommand{\redstep}{\rightarrow}
\newcommand{\redgt}{\stackrel{+}{\rightarrow}}
\newcommand{\redlt}{\stackrel{+}{\leftarrow}}
\newcommand{\redge}{\stackrel{*}{\rightarrow}}
\newcommand{\redle}{\stackrel{*}{\leftarrow}}
\newcommand{\redeq}{\stackrel{*}{\leftrightarrow}}
\newcommand{\prefixle}{\preceq}
\newcommand{\prefixlt}{\prec}
\newcommand{\prefixge}{\succeq}
\newcommand{\prefixgt}{\succ}
\newcommand{\imprefix}[1]{{\operatorname{pre}(#1)}}

\newcommand{\restr}[2]{{#1}|_{#2}}
\newcommand{\pres}[2]{\langle #1 \mid #2 \rangle}

\newcommand{\globa}{a}
\newcommand{\globb}{b}
\newcommand{\Klein}{\mathcal V_4}
\newcommand{\meet}{\wedge}

\newtheorem{mythm}{Theorem}[section]
\newtheorem{mycor}[mythm]{Corollary}
\newtheorem{mylem}[mythm]{Lemma}
\newtheorem{myprop}[mythm]{Proposition}
\newtheorem{mydefn}[mythm]{Definition}
\newtheorem{myrem}[mythm]{Remark}
\newtheorem{myexmp}[mythm]{Example}
\newtheorem{myobserve}[mythm]{Observation}

\newenvironment{rmdefn}%
        {\begin{mydefn}\rm}{\end{mydefn}}
\newenvironment{rmexmp}%
        {\begin{myexmp}\rm}{\end{myexmp}}
\newenvironment{rmrem}%
        {\begin{myrem}\rm}{\end{myrem}}
        {\begin{myprob}\rm}{\end{myprob}}
\newenvironment{rmobserve}%
        {\begin{myobserve}\rm}{\end{myobserve}}